\documentclass[10pt,journal,compsoc]{IEEEtran}

\usepackage{cite}
\usepackage{calc}
\usepackage{amsmath,amssymb,theorem}
\usepackage{graphicx}
\usepackage{mathtools}
\usepackage{cuted}
\usepackage{bbm}
\usepackage[ruled,vlined]{algorithm2e}
\usepackage{multirow}
\usepackage{epstopdf}	
\usepackage{color,xspace}
\usepackage{subfigure}
\usepackage{lipsum}
\usepackage{algpseudocode}
\usepackage{etoolbox}

\newcommand{\MarkovState}{Y}
\newcommand{\NN}{{\mathcal{N}}}
\newcommand{\until}[1]{\{1,\dots, #1\}}
\newcommand{\CT}{{x}}

\newcommand{\Si}{{g_i^{-1}}}
\newcommand{\Sj}{{g_j^{-1}}}
\newcommand{\GG}{\mathcal{G}}
\newcommand{\CTbelief}{\mathcal{\pi}}

\newcommand{\infectedneighbors}{{\mathcal{N}_i^I}}
\newcommand{\threshold}{\mathcal{T}^S}
\newcommand{\x}{\tilde{x}}
\newcommand{\C}{\mathcal{C}}
\newcommand{\R}{\mathbb{R}}
\newcommand{\infectedthreshold}{\mathcal{T}^I}
\newcommand{\V}{\mathcal{V}}
\newcommand{\U}{\mathcal{U}}

\newcommand{\oprocendsymbol}{\hbox{$\bullet$}}
\newcommand{\oprocend}{\relax\ifmmode\else\unskip\hfill\fi\oprocendsymbol}

\newcommand{\longthmtitle}[1]{\mbox{}\textup{\textbf{(#1)}}}

\newtheorem{problem}{Problem}
\newtheorem{theorem}{Theorem}
\newtheorem{remark}{Remark}

\begin{document}

\title{An Agent-Based Distributed Control of Networked SIR Epidemics}

\author{Mohammad Mubarak \qquad Cameron Nowzari

\IEEEcompsocitemizethanks{
\IEEEcompsocthanksitem The authors are with the Department of Electrical and Computer Engineering, George Mason University, Fairfax, VA 22030, USA, {\tt\small \{mmubara,cnowzari\}}@gmu.edu}}


\IEEEtitleabstractindextext{%
\begin{abstract}
This paper revisits a longstanding problem of interest concerning the distributed control of an epidemic process on human contact networks. Due to the stochastic nature and combinatorial complexity of the problem, finding optimal policies are intractable even for small networks. Even if a solution could be found efficiently enough, a potentially larger problem is such policies are notoriously brittle when confronted with small disturbances or uncooperative agents in the network. 
Unlike the vast majority of related works in this area, we circumvent the goal of directly solving the intractable and instead seek \textit{simple} control strategies to address this problem. 
More specifically, based on the \textbf{locally} available information to a particular person, how should that person make use of this information to socialize as much as possible while ensuring some desired level of safety?  
We set this up as a finite time optimal control problem using a well known exact Markov chain Susceptible-Infected-Removed (SIR) model. 
Leveraging results from the literature, we employ a commonly used mean-field approximation (MFA) technique to relax the problem. 
We show solution of this problem to be a form of threshold on the chance of infection of the neighbors of that person. Simulations illustrate our results.


\end{abstract}
}
\maketitle


\section{\textbf{Introduction}}
Modeling epidemic processes has been a longstanding research area with the earliest models proposed by Bernoulli in 1760~\cite{WOK-AGM:27,NTB:75}. Most current works in the literature focus on deterministic (generally mean-field) models that are only good for tracking aggregate numbers of infectious, while exact stochastic models (Markov Chain models) are much more suitable for understanding the spread of a disease at the person-to-person level. Although deterministic mean field models are simpler to analyze, their predictive power is poor in small networks. On the other hand, exact Markov chain model gives accurate predictions but its intractable to analyze even for small networks.  The majority of existing works only focus on the deterministic models~\cite{CN-VMP-GJP:16,FDS-CS-PVM:13,PEP-JL-CLB-BEK-TB:18,PEP-CLB-AN:17} while much fewer works that instead study the exact Markov chain compartmental models. More specifically, the work~\cite{PVM-JO-RK:08} investigates the connections between the exact Markov chain models ($2^N$ dimensional Markov Chain) and their mean-field approximations for the SIS compartmental model. In~\cite{NAR-BH:15,SO-SB:19,PVM:14} the authors extend this type of analysis to the slightly more complicated SIRS model ($3^N$ dimensional Markov chain). In~\cite{KM-MGR-CRL:16,IZK-CGM-FS-PLS-RRW-15} the exact SIR model is analyzed for various specific small graph structures or graphs with some special properties (e.g., no loops). In our previous work \cite{MM-JB-CN:21}, we studied the SEIR exact stochastic (Markov Chain) epidemic spreading model and its deterministic MFA epidemic spreading model. 
Another work that studies the the relation between stochastic and deterministic SIS, SIR in discrete time is~\cite{LJSA-AMB:00}. Many of these established results are discussed in the book~\cite[Chapter 2]{IZK-JCM-PLS:17}. 

While the above works are all interested in analyzing the connections between exact and approximate models, our ultimate goal is controlling these processes. While accurate and efficient modeling and analysis would be useful for the control problem, it is not actually necessary to establish this connection. Our goal of this paper is instead to jump directly to the intractable control problem of interest and identifying much simpler control strategies that are still effective. 
%
The reason for few works of controlling exact model is the complication that this model proposes when it comes to analyze control problems. One of the earliest works that consider a very similar setup is~\cite{KD-AO-JNT:14} where a similar exact model is considered but the available control actions came in the form of curative resources (e.g., individuals can pay some cost to recover from an infection faster such as by going to a doctor). Instead, in this paper we only consider Non-Pharmaceutical Interventions (NPIs), or more specifically the act of avoiding social interactions with chosen people.

Alternatively, most other works focus on a more direct approach of minimizing the spectral radius of an associated matrix to suppress the epidemic \cite{YW-SR-AS:07, PVM-DS-FK-CL-RB-DL-HW:11}. In \cite{AK-TB:14}, proposes an optimal control problem for a centralized network controller that regulates the infection levels in the network via adapting the curing rates of the nodes, where they used the heterogeneous SIS linearized MFA networked model. 
Optimal resource allocation problem in \cite {VMP-MZ-CE-AJ-GJP:13, CN-VMP-GJP:15}. In \cite{FL-MB:20} an optimal control problem was solved for rumor spreading on node level. In \cite{EAE-JJM-MLB:12} the authors studied an optimal link removal to minimize the spread of infection via quarantining with limited resources. Other work proposed algorithm with approximation, to to minimize the number of infected people, based on the idea of bounding the number of infection by supermodular function. \cite{YY-LS-PEP-KHJ:20}.  An optimal control problem with data-driven model for the spread of COVID-19  and minimize the economic costs associated with implementing NPIs \cite{MH-FB-VMP:20}.
In \cite{YM-SC:13}, the authors proposed an optimal control problem formulation to minimize the total number of infectious during the spread of SIR epidemics by controlling the contact rate, they ended up solving the problem for centralized well-mixed homogeneous SIR model. 
In \cite{CC-ZL-CM-CM:21} a centralized optimal control problem to control the contact rate between agents was presented and solved numerically using SQP. More related to our problem, the  authors in \cite{ER-HE-TG:20} studied an optimal control problem for an individual who is trying to avoid getting infected by controlling their contact rate. However, they considered a well mixed homogeneous SIR model and solved the problem numerically. Another work proposed in \cite{SYO-SA-JWK-EM-TB-MW-PGM:22} where an agent controls its rate of contact with others such that partial observability of viral status is considered. 

%

\textit{Statement of Contributions:} First, we formulate a distributed optimal control problem, where we approach the problem using bottom-up strategy by formalizing a local stochastic optimal control problem from the point of view of a single person in a social network. We then consider the mean-field approximation of the problem for which a control strategy is proposed where interactions with certain people are stopped when their chance of infection exceeds some threshold. We show that the optimal solution to the relaxed problem is a controller of this form. We then verify the effectiveness of our solutions on the original stochastic problem, rather than the relaxed problem. 
Finally, we take the problem a step further, by analyzing the solution of a single person interacting with a lumped population, we found that the optimal control strategy also depends on the chance of infection of the lumped node, such that, the single individual will fully quarantine if the chance of infection of the lumped node exceed some threshold. We show numerically how such a strong assumption can increase the cost of disconnection on the individual from considering a heterogeneous population. 

\textit{Preliminaries and notations:} We denote by $\R$, $\R_{>0}$, $\R_{\ge 0}$, and $\mathbb{Z}_{\ge0}$ the set of real, positive real, non-negative real, and non-negative integer numbers, respectively. 

We denote an n-dimensional column vector with each entry equal to $1$ by $\mathbf{1}_{n}$. 
The cardinality of a set $v \in \R^n$ is denoted by $|v|$. We say that the matrix $A\in \R^{n\times n}$ is symmetric if $A=A^T$. If a vector $x \in \R^N$, we denote the diagonal matrix of $x$ by diag$(x)\in \R^{N\times N}$, where all the off diagonal entries are zeros while the main diagonal contains all the elements of $x$. The Cartesian product of two sets $A$ and $B$ is denoted by $A\times B=\{(a,b)|a\in A, b \in B\}$ which represents the set of all points $(a,b)$, where $a\in A$ and $b\in B$.





\quad \textit{Graph Theory:} 
an unweighted undirected graph $\GG = (\mathcal{V}, \mathcal{E}, A)$, where the set of vertices~$\mathcal{V} = \until{N}$ 
captures all the nodes in a network, the edge set~$\mathcal{E} \subset \mathcal{V} \times \mathcal{V}$ denoting the interactions between the different nodes, and the adjacency matrix~$A = [a_{ij}] \in \R^{N \times N}$, where~$a_{ij} = 1$ if $(i,j) \in \mathcal{E}$, 
and $a_{ij} = 0$, otherwise. Unidrected graphs implies a symmetric adjacency matrix $A$. 

\section{\textbf{Problem Formulation}}


Consider a group of $N$ people (i.e., nodes), interacting according to an unweighted undirected graph $\GG=(\mathcal{V}, \mathcal{E}, A)$. An interaction, in this problem context, is any activity that reguarly brings individuals close enough to spread a disease. If $(i,j)\in \mathcal{E}$, then person~$i$ and person~$j$ have contact with each other and the disease can spread between them. Imagining the problem from the point of view of the whole network, we are interested in analyzing and formulating a stochastic optimal control problem trading off the people's desire to socialize against the desire to stay healthy. 
More specifically, each person should socialize as much as possible while ensuring that their probability of remaining healthy after some period of time~$T$ is greater than some personal threshold~$\threshold_1 \in [0,1]$. A more cautious person would select a larger~$\threshold_1$ for instance.

\subsection{Compartmental Modeling}\label{se:modeling}
\begin{figure}[h]
	\centering
	\includegraphics[width=0.9\linewidth]{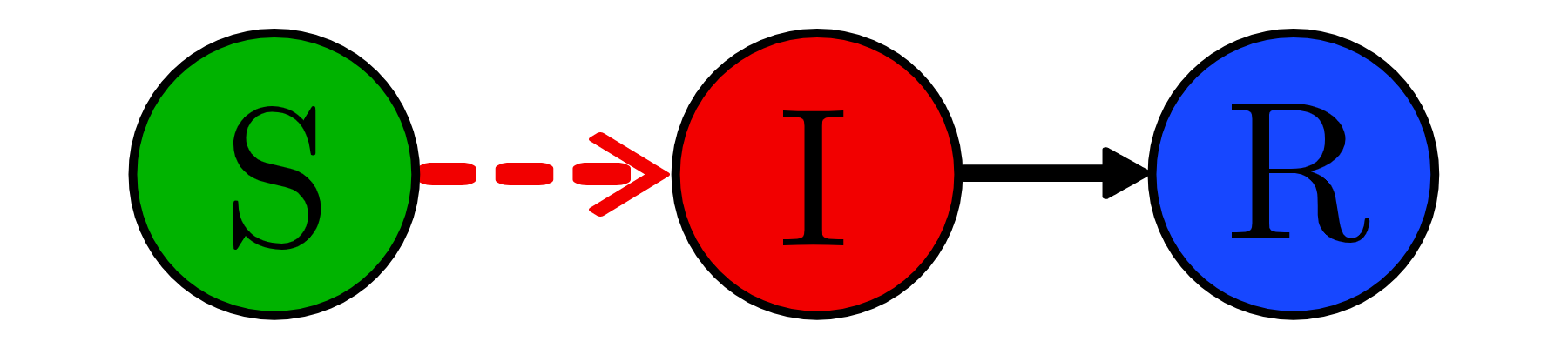}
	\put(-80,31){$\delta$}
	\put(-160,34){$\beta^\text{eff}_i(t)$}
	\caption{Graphical representation of SIR compartmental model for a single person. }
\label{fig:transition}
\end{figure}

We use a Markov Process to keep track of the exact Markov states and the state of each node $i \in \mathcal{V}$ in the whole network.
We start with Figure~\ref{fig:transition} which shows the Susceptible-Infected-Removed (SIR) compartmental model for a single person. A person in the Infected state will naturally move to the Removed state over time. However, a person that is in the Susceptible compartment can only move to the Infected compartment through interactions with infected individuals. 
 
The rate at which an individual might transition from one compartment to the next are defined by~$\beta^\text{eff}_i(t), \delta > 0$, where~$\beta^\text{eff}_i(t)$ depends on the states of the individual's neighbors. The term~$\beta^\text{eff}_i(t)$ will be explained soon, whereas the recovery rate~$\delta$ is a fixed constant that doesn't depend on interactions with other people. 
Note that our model does not distinguish between people who have recovered or have died and we lump these individuals in the `Removed' state. 
The Markov process for the entire system is denoted by $[X(t)]_{t\ge0}$, where $X(t)=[X_1(t),\dots, X_N(t)] \in \C^N\triangleq \{S, I, R \}^N$ is the entire Markov state and~$X_i(t) \in \C$ denotes the state of individual~$i$. 
The total number of possible states of the network is then $n=3^N$.

Now, for convenience, we define $\NN_i \triangleq \{j \in \mathcal{V}|(i,j)\in \mathcal{E}\}$ as the set of all neighbors of node $i$. We also define $\infectedneighbors \triangleq \{j \in \mathcal{N}_i|X_j(t) = I\}$ 
as the set of infected neighbors of node $i$.
The total infection rate felt by node~$i$ is then given by $\beta^\text{eff}_i(t) \triangleq \beta \sum_{j \in \infectedneighbors} (1-\operatorname{max}(u_{ij}(t),u_{ji}(t)))$, where $u_{ij}(t) \in [0,1]$ is the control input which will be formalized in the next section. In other words, each infected neighbor contributes~$\beta > 0$ to the total rate of infection when they are interacting normally. 

The Markov process is defined by the following Poisson rates, as presented graphically in Figure~\ref{fig:transition}, 
\begin{align*}
&X_i:\text{S}\rightarrow \text{I with rate } \beta^\text{eff}_i(t),\\
&X_i:\text{I}\rightarrow \text{R with rate } \delta.
\end{align*}

\subsection{Control Mechanism}\label{se:control}

We formulate the real world problem for the entire network as a distributed optimal control. In order for each node~$i$ to be able to protect themselves, we allow the ability to disconnect their own links from neighbors (e.g., avoiding contacts with specific people); or even disconnecting entirely (e.g., quarantine) when needed. To model this, we define the weighted undirected subgraph~$\hat{\GG}(U(t))=(\mathcal{V},E(t),\hat{A}(t))$,
which  has the same set of vertices of the static graph $\GG$, but a different set of edges $E(t) \subset \mathcal{E}$, and weighted adjacency matrix $\hat{A}(t)=(A-U(t))\in [0,1]^{N \times N}$, where $U(t)=[u_{ij}(t)] \in [0,1]^{N \times N}$ is defined as the action matrix, where ${u}_{ij}(t)=0$ if $a_{ij} = 0$.  If~$a_{ij} = 1$, then the value~$u_{ij} \in [0,1]$ represents the level of caution between persons~$i$ and~$j$, where the link $(i,j)$ will depend on the more cautious person, with~$u_{ij} = 1$ meaning that they are not interacting at all and thus the infection cannot spread between the two. Further, we define the set $\U_i=\{u_{ij}\in[0,1]|j\in \NN_i\} \in [0,1]^{|\NN_i|}$,
which captures all the control inputs of node ``i'', and the set $\U=\U_1 \times \U_2 \times \dots \times \U_N$, which captures all the control inputs for the whole network.

{
The global cost, (e.g., global social cost)  of using this control is given as a linear function of this input so the total cost of interest to minimize is
\begin{align}\label{eq:globalcost}
J_G = \int_0^T \sum_{i\in \V}\sum_{j \in \NN_i} c_{ij} \operatorname{max}(u_{ij},u_{ji}) dt 
\end{align}
where~$T > 0$ is the time horizon, and $c_{ij}=c_{ji} > 0$ is the cost associated with disconnecting edge $(i,j)$ from the graph.

The question then is how should each person~$i$ actively control their own links~$u_{ij}(t)$ for~$j \in \NN_i$ to minimize the global cost~\eqref{eq:globalcost} while satisfying the constraint that its probability of being healthy at any time~$t\in[0,T]$, (i.e., over some time horizon), is greater than or equal to~$\threshold_i$,
\begin{align}\label{eq:pisafety}
\operatorname{Pr}[X_i(t) = S] \geq \threshold_i.
\end{align}


This is formalized in Problem~\ref{pr:globalreal-world problem}. 

\begin{problem}[Distributed Control Problem]\label{pr:globalreal-world problem}
{\rm
\begin{align*}
\underset{u_{ij}\in \U}{\text{minimize}}\quad & \int_0^T  \sum_{i\in \V}\sum_{j \in \NN_i} c_{ij} \operatorname{max}(u_{ij},u_{ji}) dt ,\nonumber\\
\text{subject to}& \quad \operatorname{Pr}[X_i(t) = S] \geq \threshold_i,\, \forall i \in \V,\nonumber\\
& \quad \text{ and dynamics}.
\end{align*} 
}
\end{problem}
To this end, describing the dynamics constraint can change the complexity of solving Problem~\ref{pr:globalreal-world problem}. One may assume a well mixed homogeneous network (all interacting with all), which can greatly simplify the problem. However, since we are dealing with person-to-person level, well mixed assumption is not very realistic, thus we avoid such an assumption in our formulation (except in Section~\ref{se:well-mixed}, to show why such an assumption in not good). Thus, since we count the heterogeneous interactions in our formulation, solving Problem~\ref{pr:globalreal-world problem} using top-down approach, (i.e., approaching the macroscopic problem as whole), is complicated and can have scalability issues. 
That said, we need to tackle the complication and scalability issue of this distributed problem in a different way. We propose a bottom up approach (i.e., identify and solve microscopic pieces of the macroscopic problem), to solve Problem~\ref{pr:globalreal-world problem}. Such that we formalize a local optimal control problem for each individual in the network, based on their own personal needs. As a result, the combination of these local strategy will impact the global quantity. Indeed, local, from individual perspective, means that each individual cares more about their selves and their objective (e.g., social interaction) is more of something that can be beneficial for them.  
}

To this end, we start to formulate the problem from the point of view of each node, e.g., node~$1$(without loss of generality). Since node~1 is the one trying to satisfy its safety, thus, we consider the control input related to node~1, (i.e., we set $u_{ij}=0$ for all $i \in \{2,\dots,N\}$). The social cost of using this control is given as a linear function of this input so the total cost of interest to minimize is


\begin{align*}
J = \int_0^T \sum_{j \in \NN_1} c_{1j}\operatorname{max}(u_{1j},u_{j1}) dt 
\end{align*}
where~$T > 0$ is the time horizon, and $c_{1j} > 0$ is the cost associated with disconnecting edge $(1,j)$ from the graph. 
Since we've set $u_{ij}=0$ for all $i \in \{2,\dots,N\}$, 
the social cost function becomes,

\begin{align}\label{eq:cost}
J = \int_0^T \sum_{j \in \NN_1} c_{1j} u_{1j} dt 
\end{align}

The question then is how should person~$1$ actively control their own links~$u_{1j}(t)$ for~$j \in \NN_1$ to minimize the cost~\eqref{eq:cost} while satisfying the constraint that its probability of being healthy at any time~$t\in[0,T]$, (i.e., over some time horizon), is greater than or equal to~$\threshold_1$,
\begin{align}\label{eq:p1safety}
\operatorname{Pr}[X_1(t) = S] \geq \threshold_1.
\end{align}

This is formalized in our main Problem~\ref{pr:real-world problem}. 

\begin{problem}[Maintaining Individual Safety]\label{pr:real-world problem}
{\rm
\begin{align*}
\underset{u_{1j}\in \U_1}{\text{minimize}}\quad & \int_0^T \sum_{j \in \NN_1} c_{1j}u_{1j}(t) dt ,\nonumber\\
\text{subject to}& \quad \operatorname{Pr}[X_1(t) = S] \geq \threshold_1,\nonumber\\
& \quad \text{ and dynamics}.
\end{align*} 
}
\end{problem}

%

\begin{table}
\begin{center}
\begin{tabular}{|cl|}
\hline
variable & description \\ 
\hline
$N$ & number of individuals in the network \\
$n = 3^N$ & number of discrete Markov states \\
$c_{ij} \in \R_{\geq 0}$ & cost of reducing contact with agent~$j$ \\
$\beta_i^\text{eff}(t)>0$ & total rate of infection of agent $i$ \\
$\beta>0$ & transmission rate \\
$\delta>0$ & recovery rate \\
$X \in \C^N$ & Markov state  \\
$X_i \in \C$ & compartmental state of agent $i$ in Markov state $X$ \\
$\CT_i^S $ & exact probability of agent $i$ being Susceptible \\
$\x_i^S $ & relaxed MFA probability of agent $i$ being Susceptible\\
$u_{ij} \in [0,1]$ & control input between agent $i$ and agent $j$\\
$\U_i \in [0,1]^{|\NN_i|}$& the set of all control inputs of node ``i''\\
$\threshold_i \in [0,1]$ & personal threshold of agent $i$\\
$\infectedthreshold_i \in [0,1]$ & infected threshold of agent $i$\\
$\NN_i$ & set of all neighbors of agent $i$\\
$\infectedneighbors$ & set of infected neighbors of agent $i$\\
\hline
\end{tabular}
\end{center}
\caption{Agent (node)~$i$ model definitions}\label{ta:glossary}
\end{table}


\section{\textbf{Optimal Control Formulation}}\label{se:Main_Problem_Sol}

Here we expand out the dynamic constraints of Problem~\ref{pr:real-world problem} to formalize the mathematical Problem~\ref{pr:SC}. 

\subsection{Continuous Time Exact Markov Chain Model (CT-Exact)}\label{se:CTMP}
Since we are considering the spread at individual level we formalize the networked Continuous-Time Exact Markov Chain/Process (i.e., stochastic) SIR model of interest. Indeed, this model is not new and it's known in the literature, but is notoriously intractable to analyze~\cite{NAR-BH:15,SO-SB:19,PVM:14,PVM-JO-RK:08, MY-CS:11}.
Now, in order to analyze the Markov process exactly, we need to consider the probability distribution across the total number of possible states in the entire network~$n = 3^N$. We refer to the collection of all distinct states as state space denoted $\mathcal{Y} =\{\MarkovState_1,..,\MarkovState_n\}$, 
where each Markov state $\MarkovState \in \mathcal{Y}$ corresponds to exactly one element in~$\C^N$.
Define $g:\C^N \mapsto \mathcal{Y}$ as the bijective map from each state~$X(t) \in \C^N$ to a unique Markov state~$\MarkovState(t) = g(X(t)) \in \mathcal{Y}$. Then given a Markov state~$\MarkovState \in \mathcal{Y}$, we can extract the compartmental state of node~$i$ by~$g_i^{-1}(\MarkovState) \in \C$.

Since this is a probabilistic model, let $\CTbelief(t) \in [0,1]^{n}$ be the belief state (i.e., probability distribution) of the Markov process, capturing the probability of the Markov process being in each of the states, i.e.,~$\CTbelief_m(t) = \text{Pr}[g(X(t))=\MarkovState_m]$.
Since this is a probability distribution, we have the property $\mathbf{1}_n^T\CTbelief(t)=1$ at all times. 

Now, given some initial belief state $\CTbelief(0)=\CTbelief_0$, we wish to propagate the belief state forward in time, so that we know what the probability of being in each state is. This can be accomplished using the transition matrix
of the time varying graph $P(U(t)) \in \R^{n\times n}$ of the continuous time Markov process, where only one node $i\in \mathcal{V}$ can change state at a time. We write $P$, as a function of $U(t)$, to emphasis that it explicitly depends on the control input. Each entry $p_{rc}$, in the transition matrix, represents the rate of change from Markov state $Y_r$ to Markov state $Y_c$ is defined by 

\begin{align*}
p_{rc}=
&\begin{cases}
      \delta, &  \text{if}\ \Si(\MarkovState_{r})=I,\, \Si(\MarkovState_{c})=R,\nonumber\\&\Sj(\MarkovState_{r})=\Sj(\MarkovState_{c})$, for $j\neq i. \nonumber \\
      \beta^{\text{eff}}_i(t),&\text{if}\ \Si(\MarkovState_{r})=S,\, \Si(\MarkovState_{c})=I,\nonumber\\& \Sj(\MarkovState_{r})=\Sj(\MarkovState_{c})$, for $j\neq i.\nonumber \\ 
      -\sum_{b=1:b\neq c}^n p_{bc},&\text{if }r=c.\nonumber \\
      0,&\text{otherwise},
\end{cases}
\end{align*}

where $\Si(\MarkovState_{c})$, $\Si(\MarkovState_{r})$ are the compartmental state of node $i$ in Markov state $\MarkovState_c$ and $\MarkovState_r$, respectively , for any node $i \in \mathcal{V}$, and $c,r=1,\dots,n$. 

Now, we start writing the forward propagation equation of continuous time exact Markov chain model by applying the law of total probability,

\begin{align}\label{eq:limit}
\dot \CTbelief_c(t) &=\sum_{r=1}^n {\CTbelief_r(t){p_{rc}}},
\end{align}
Now, we can define the forward propagation of the belief state $\CTbelief(t)$ in matrix form as
\begin{align}\label{eq:3 state MC}
\dot \CTbelief^T(t)&=\CTbelief^T(t)P(U(t)).
\end{align}


In order to extract the probability of a given node~$i$ being in a compartmental state~$S, I, R \in \C$, we simply sum over the distribution of Markov states~$Y \in \mathcal{Y}$ that represent~$X_i = S$, $X_i = I$, or $X_i = R$, 
\begin{align}\label{eq:marginalProbs}
\CT_i^S(t) \triangleq \operatorname{Pr}[X_i(t) = S] = \sum_{m \in \until{n} | \Si(Y_m) = S} \pi_m,\\
\CT_i^I(t) \triangleq \operatorname{Pr}[X_i(t) = I] = \sum_{m \in \until{n} | \Si(Y_m) = I} \pi_m, \\
\CT_i^R(t) \triangleq \operatorname{Pr}[X_i(t) = R] = \sum_{m \in \until{n} | \Si(Y_m) = R} \pi_m .
\end{align}
Where $\CT^S_i(t)+\CT^I_i(t)+\CT^R_i(t)=1$ holds at all times, due to the probabilistic nature of the model.

We can now formalize the constraint of Problem~\ref{pr:real-world problem} and present this in~Problem~\ref{pr:SC}. 

\begin{problem}[Agent-Based Stochastic Optimal Control]\label{pr:SC}
{\rm
\begin{align*}
\underset{u_{1j}\in \U_1}{\text{minimize}}\quad & \int_0^T \sum_{j \in \NN_1} c_{1j}u_{1j}(t) dt ,\nonumber\\
\text{subject to} \quad & \CT_1^S(t) \geq \threshold_1,\nonumber\\
&\dot \CTbelief^T(t)=\CTbelief^T(t)P(U(t)),\nonumber\\
&\CTbelief(0)=\CTbelief_0,\nonumber\\
& u_{1j}(t)\in [0,1], \quad \forall j \in \NN_1, t\in[0,T].
\end{align*} 
}
\end{problem}

It is worth noting that the problem is feasible as long as the probability of node~$1$ being Susceptible is initially greater than or equal to the threshold~$\threshold_1$. This is easy to see as simply disconnecting entirely from all neighbors~$u_{1j}(t) = 1$ will satisfy the constraint (but with a worst-case cost). 

%

%
%

\begin{remark}[POMDP]
{\rm
Problem \ref{pr:SC} can also be looked at as a Partially Observable Markov Decision Process (POMDP), with the tuple of $(\mathcal{Y}, \U, P(U(t)), J,\mathcal{O}, B)$ where $\mathcal{Y}$ is the state space, $\U$ is the action (control) set, $P(U(t)$ is the transition matrix, $J$ is the decision running cost, which were predefined. While $\mathcal{O}$ represents a finite observation set, and $B$ observation function, which are not considered in this work scoop. General talking, Covid-19 Antigens testing results can be considered as the observation set for all people on network. This formulation and its solution techniques is beyond this work scoop. But it's worth to be mentioned since it can be similar to the problem formulation~\ref{pr:SC}. For more details about POMDP we refer reader to \cite[Part-II]{VK:16}. \oprocend
}
\end{remark}

We acknowledge that when dealing with a person level, the predictions using Markov chain exact model is more accurate to use than MFA \cite{MM-JB-CN:21}. However, solving the control Problem~\ref{pr:SC} exactly is generally intractable and hard to analyze due to the~$O(3^N)$ size complexity of the problem, but is indeed of great interest. To the best of our knowledge, no one have ever solved the original problem~\ref{pr:SC} before or even show the effectiveness of the solution to the relaxed problem on the original problem. Instead for now, we resort to relaxing the problem to provide a sub-optimal solution to Problems~\ref{pr:real-world problem} and~\ref{pr:SC} instead.

\section{\textbf{Sub-optimal Solution Approach}}
In order to obtain a sub-optimal feasible solution to Problem~\ref{pr:SC}, we first relax the problem and then solve the relaxed problem optimally by proposing a simplifying assumption. While this is not as good of a solution we hope to eventually obtain, we believe it is the best solution that currently exists in the literature for Problem~\ref{pr:SC} as the majority of works consider relaxed problems to begin with (e.g., mean-field or lumped degree models). To leverage existing results, we similarly perform a Mean Field Approximation to relax the problem~\cite{PVM:14,PVM:11}, but unlike other works we evaluate the effectiveness of its solution on the original Problem~\ref{pr:SC}.

\subsection{Continuous Time Mean Field Approximation (CT-MFA)}\label{se:CTMFA}

The MFA approximation aims to reduce the number of states from $3^N$ to $2N$, by using a first order moment closure technique. In other words, it aims to reduce the size complexity of the problem form being exponential $O(3^N)$ to linear $O(2N)$. 
This can be done by assuming an independent random variables among the nodes in a network. In order to derive the MFA for the exact model, we need to close the exact model \eqref{eq:marginalProbs}, where we consider the first moment closure, which results the first order MFA of the continuous time exact model. This relaxation is known as the ``N-Intertwined Mean field approximation (NIMFA)'' \cite{PVM:11}. 
Note that we use ( \textbf{$\tilde{}$} ) notation to differentiate between the relaxed and the exact states.
\begin{align}\label{eq:CTMFA}
\dot{\x}_i^S&=-\x_i^S(t)\sum_{j\in\NN_i}{\beta\left(1-\operatorname{max}(u_{ij}(t),u_{ji}(t))\right)}\x_j^I(t),\nonumber\\
\dot{\x}_i^I&=-\delta {\x}_i^I(t)+{\tilde{x}}_i^S(t)\sum_{j\in\NN_i} {\beta\left(1-\operatorname{max}(u_{ij}(t),u_{ji}(t))\right)} {\x}_j^I(t),\nonumber\\
\dot{\x}_i^R&=\delta {\x}_i^I(t).
\end{align}

Due to the probabilistic nature of the model, $\x_i^S(t)+\x_i^I(t)+\x_i^R(t)=1$ will always be valid for all time $t$. Thus, one can reduce the number of the states to $2N$, such that $\x_i^R(t)=1-\x_i^S(t)-\x_i^I(t)$. This drops the need of keep tracking the recovered states. 

%
\begin{theorem}\longthmtitle{Relation between exact~\eqref{eq:marginalProbs} and relaxed (MFA) models~\eqref{eq:CTMFA}~\cite{EC-PVM:14,RCB-LB-PLS-IZK:19}}\label{th:Inequality}

At any time $t \in [0,T]$, the solution, to the relaxed model~\eqref{eq:CTMFA}, $\x_1^S(t)$, will always satisfy the personal threshold in Problem \ref{pr:SC}. Such that
\begin{align}\label{eq:inequality}
\CT_1^S(t) &\ge {\x}_1^S(t).
\end{align}
\end{theorem}

\begin{IEEEproof}
It's known that for any $i\in \mathcal{V}$ and for all time $t \in [0,T]$, the probability of being Susceptible, ${\x}_i^S(t)$, in the relaxed MFA model \eqref{eq:CTMFA} lower bounds the exact one, $\CT_i^S(t)$, in \eqref{eq:marginalProbs}, \cite{EC-PVM:14,RCB-LB-PLS-IZK:19}. 
Which can be written mathematically as 
\begin{align}
\CT_i^S(t) &\ge {\x}_i^S(t),\label{eq:bnd1}
\end{align}
This leads to the validity of the inequality~\eqref{eq:inequality}.
\end{IEEEproof}

Theorem~\ref{th:Inequality} validates the use of MFA model to describe the dynamics. Even though, the predictions of MFA is not as accurate as exact model, thanks to this inequality which guarantees the safety threshold $\CT_1^S(t) \ge {\x}_1^S(t) \ge \threshold_1$.  
Another note here, as the number of agents $N$ on a network increases, the solution to MFA relaxed model~\eqref{eq:CTMFA} should become closer to the exact solution of \eqref{eq:marginalProbs}, such that as $N\rightarrow \infty$ the exact solution should match the relaxed MFA solution. This gives a good reason of the popularity of the MFA models over the exact model, which is intractable and computationally costly for a network that contains large number of agents.
For more details about the relation between the exact and MFA models we refer the reader to our previous work in \cite{MM-JB-CN:21}. 

Now, we write the relaxed problem, from point of view of node ``1'', as follows,

\begin{problem}\label{pr:CT-MFA}\longthmtitle{Relaxed Agent-Based Stochastic Optimal Control}
\end{problem}
\begin{align}\label{eq:Opt-CTMFA}
\underset{u_{1j}\in \U_1}{\text{minimize}}\quad & \int_0^T \sum_{j \in \NN_1} c_{1j}u_{1j}(t) dt ,\nonumber\\
\text{subject to} \quad & \x_1^S(t) \geq \threshold_1,\nonumber\\
&\dot{\x}_i^S=-\x_i^S(t)\sum_{j\in\NN_i}{\beta\left(1-u_{ij}(t)\right)}\x_j^I(t),\nonumber\\
&\dot{\x}_i^I={\tilde{x}}_i^S(t)\sum_{j\in\NN_i} {\beta\left(1-u_{ij}(t)\right)} {\x}_j^I(t)-\delta {\x}_i^I(t), \nonumber\\
&\x(0)=x_{0}, \, i \in \mathcal{V}, \nonumber\\
& u_{1j}(t)\in [0,1], \quad \forall j \in \NN_1, t\in[0,T].
\end{align}
where $\x(0)=[[\x^S_1(0),\dots, \x^S_N(0)]^T, [\x^I_1(0),\dots, \x^I_N(0)]^T] \in \R^{N \times 2}$. 
This initialize all the nodes in the network. 

An important note here, based on the results in Theorem~\eqref{th:Inequality}, the optimal solution to the relaxed Problem~\eqref{pr:CT-MFA} will solve the the original Problem~\eqref{pr:SC} sub-optimally.

\begin{theorem}\longthmtitle{Problem \ref{pr:CT-MFA} solution form}\label{th:solut1}
The form of optimal solution of Problem \ref{pr:CT-MFA}, is bang-bang controller, where the input form
\begin{align}\label{eq:switchingfun} 
u_{1j}(t)=&
\begin{cases}
1,\quad &\text{if } c_{1j}+\sum_{i \in \mathcal{V}}(\lambda_{\x_i^{S*}}-\lambda_{\x_i^{I*}}) \x_i^{S*} \x_j^{I*} < 0,\\
0, \quad &\text{if } c_{1j}+\sum_{i \in \mathcal{V}}(\lambda_{\x_i^{S*}}-\lambda_{\x_i^{I*}}) \x_i^{S*} \x_j^{I*}>0.
\end{cases},
\end{align}

\end{theorem}

\begin{IEEEproof}
See appendix.
\end{IEEEproof}

Note that in order to find closed form solution of the optimal input $u_{1j}^*$, we need to get a closed form solution of the states $\x_i^S, \, \x_i^I$ for all $i \in \V$, and the co-states $\lambda$, by methods of solving odes. Getting such a closed form solutions are generally intractable \cite{AK-TB:14}. However, one can turn to numerical optimization techniques to get an approximated numerical solution to problem~\eqref{pr:CT-MFA}. The problem~\eqref{pr:CT-MFA} falls under the category of non-linear constrained optimization problems, which is known to be very hard to solve. 

To this end, we solve problem~\eqref{pr:CT-MFA} in two ways. First, we use OpenOCL since it provides a modeling language that helps in solving constrained optimal control problems numerically. Open Optimal Control Library (OpenOCL)\cite{koenemann2017openocl}, is a software package that implements direct collocations methods and Ipopt, and interfaces CasADi-software tools \cite{JAEA-JG-GH-JBR-MD:19} to solve a non-linear program.  
Second, we present the main result in the paper, where we show an analytical form of the solution to problem~\eqref{pr:CT-MFA}. Although it's not a closed form solution, but is indeed of great interest, where it allows solving the problem in different way.


\subsection{One-way Infection to Node 1}

\begin{remark}[Simplifying Assumption]
{\rm
This assumption is motivated by the fact that, for the problem of keeping node 1 safe, we do not wish to consider cases where infection spreads from node 1, because then the safety of that node would already have been compromised. \oprocend
}
\end{remark}

Using the simplifying assumption that node~$1$ cannot infect its neighbors, but can be infected by them, we have the dynamics
\begin{align}
\dot{\x}_1^S&=-\x_1^S(t)\sum_{j\in\NN_1}{\beta}(1-u_{1j}(t))\x_j^I(t),\label{eq:dynamics1}\\
\dot{\x}_i^I&=-\delta \x_i^I(t)+\x_i^S(t)\sum_{j\in\NN_i\setminus 1}{\beta}\x_j^I(t),\label{eq:dynamics2}\\
\dot{\x}_i^S&=-\x_i^S(t)\sum_{j\in\NN_i\setminus 1}{\beta}\x_j^I(t)\label{eq:dynamics3},
\end{align}
where $i \neq 1$, such that  $i \in \{2,\dots,N\}$. 

This greatly simplifies the problem, because~\eqref{eq:dynamics2} and~\eqref{eq:dynamics3} do not depend on $\x_1^S$ or $u_1$. 
Note that one would need to consider a digraph to model this one way infection. 
Also note that the dynamics size will be at most $2N-1$ for the complete graph, since we don't need to keep track of the recovered state of any node nor the infected state of node 1. 

\begin{theorem}\longthmtitle{Problem \ref{pr:CT-MFA} solution}\label{th:solut}
\rm
{
We can solve Problem \ref{pr:CT-MFA}, with the dynamics constraints replaced by~\eqref{eq:dynamics1},~\eqref{eq:dynamics2},~\eqref{eq:dynamics3}, using the input form 
\begin{align}\label{eq:heu} 
u_{1j}(t)=&
\begin{cases}
1,\quad &\text{if } \x_j^I(t)>\infectedthreshold_j,\\
0, \quad &\text{if } \x_j^I(t)<\infectedthreshold_j.
\end{cases},
\end{align}
where $\mathcal{T}_j^I \triangleq \frac{-c_{1j}}{\beta \threshold_1 v}>0$ is a threshold on the probability of node 1's neighbor $j$ being infected, for $j \in \NN_1$, and $v \in \R$.
}
\end{theorem}
\begin{IEEEproof}
See appendix.
\end{IEEEproof}

\begin{remark}[Optimal Infection Threshold]
{\rm
Note that in this work we show that such a threshold exist, but we do not yet have a a closed form solution for the optimal threshold $\infectedthreshold_j$ for which~\eqref{eq:heu} is optimal, thus we rely on numerical optimization techniques to find this threshold $\infectedthreshold_j$. \oprocend
}
\end{remark}

To this end, to find the threshold, $\infectedthreshold_j$,  numerically we apply MonteCarlo method. Where we iterate on randomly sampled threshold $\infectedthreshold_j$ until we find this threshold. In simulation section, we will show a comparison between the time complexity of our MonteCarlo method against OCL tool.


\section{Well-Mixed Homogeneous Network}\label{se:well-mixed}
As we presented in the literature review, the majority of the researches focus the work on the well-mixed homogeneity assumption. That is for a good reason, it can greatly simplifies the number of states to track, and can solve the scalability issues for large network. Therefore, we would like here to propose an assumption that our network is well mixed, to see how bad such an assumption can be compared to the networked relaxed MFA model. 

\begin{remark}[One Control Input]
{\rm
We propose here a new constraint that node``1'' interactions can be controlled by controlling one input rather than the multi-inputs. This can be achieved by assuming a well mixed homogeneous network, such that node``1'' is interacting with a lumped population. \oprocend
}
\end{remark}

Note that, it's convenient here to consider a well mixed homogeneous network. The assumption of a well mixed homogeneous network, considers that every individual in the network is in contact with all other individuals, i.e., complete graph. This means, that node ``1'' is connected to all the individuals on the network. The well-mixed homogeneous network model is popular in the literature due to its simplicity, some of the optimal control problems that considers this model can be found in \cite{YM-SC:13, ER-HE-TG:20}. Such assumption can greatly reduce the computation cost of finding multi-input's optimal solution, specifically for large networks. However, simulations will show such an assumption is not as good for cost.

Now, we denote the lumped node's (population) average probabilities of being susceptible, infected, and recovered by $S_{LP}(t)$, $I_{LP}(t)$, and $R_{LP}(t)$, respectively. Due to the complete graph assumption in the lumped node, we can write the equations that describes the states of the lumped node as follows
\begin{align*}
S_{LP}(t) &= \sum_{i \in \mathcal{V} \setminus 1} \frac{\x_i^S(t)}{N-1},\\
I_{LP}(t) &= \sum_{i \in \mathcal{V} \setminus 1} \frac{\x_i^I(t)}{N-1},\\
R_{LP}(t) &= \sum_{i \in \mathcal{V} \setminus 1} \frac{\x_i^R(t)}{N-1}.
\end{align*}
Since we are dealing with average probabilities, thus, $S_{LP}(t)+I_{LP}(t)+R_{LP}(t)=1$ holds at any time $t$, due to the probabilistic nature of the states in the lumped node. This can justify dropping the need of keep tracking of the average probability of recovered individuals in the lumped node.
Using this simplifying assumption, i.e., homogeneously mixed network, we can write the dynamics of the system we wish to control,

\begin{align}\label{eq:node-vs-LN}
\dot{\x}_1^S&=-\beta \x_1^S (1-u) I_{LP}, \nonumber\\
\dot{I}_{LP}&=-\delta I_{LN}+\beta S_{LP} (I_{LP}+ (1-u)\x_1^I), \nonumber\\
\dot{S}_{LP}&= -\beta S_{LP} (I_{LP}+ (1-u)\x_1^I), \nonumber \\
\dot{\x}_1^I&=-\delta \x_1^I+\beta (1-u)\x_1^S I_{LP}.
\end{align}
Note that due to homogeneously mixed network assumption, the index of the control $u$ is dropped, since only one input is being controlled. 

We notify here, the problem~\eqref{pr:CT-MFA} using the constraint dynamics~\eqref{eq:node-vs-LN} can also be solved numerically using Open Optimal Control Library(OpenOCL). 

Note that we can get the form of the optimal control solution if we assume here that node``1'', is safe, thus, the chance node ``1'' spreading the infection to the lumped node can be neglected since it won't have a significant effect in spreading the infection to the lumped node. We also can drop the index on the control input since we only have one input in this problem. Considering this assumption we can simplify the dynamics in \eqref{eq:node-vs-LN} as follows, 

\begin{align}\label{eq:simplyfied-node-vs-LN}
\dot{\x}_1^S&=-\beta(1-u) \x_1^S I_{LP}, \nonumber\\
\dot{I}_{LP}&=-\delta I_{LP}+\beta S_{LP} I_{LP}, \nonumber\\
\dot{S}_{LP}&= -\beta S_{LP} I_{LP}.
\end{align}

Also the cost function becomes  
\begin{align*}
J = \int_0^T c_{av}u(t) dt,
\end{align*}
where $c_{av}=\frac{N-1}{|\NN_1|}\sum_{j\in \NN_1}c_{1j}$, is the average of all edge's associated cost. Note that in $c_{av}$, we multiply by $N-1$, due to complete graph assumption. Means that node ``1'' will be considered to have a link that connects with all neighbors.

\begin{theorem}\longthmtitle{Problem \ref{pr:CT-MFA} solution with strong assumption}\label{th:solution2}
\rm
{
We can solve Problem \ref{pr:CT-MFA}, with the dynamics constraints replaced by~\eqref{eq:simplyfied-node-vs-LN} using the input form 
\begin{align}\label{eq:heuLP} 
u^*(t)=&
\begin{cases}
1,\quad \text{if } I_{LP}(t)>I^*,\\
0, \quad \text{if } I_{LP}(t)<I^*.
\end{cases},
\end{align}
where $I^* \triangleq \frac{-c_{av}}{\beta \threshold_1 v}>0$ is a threshold on the probability of node 1's lumped node neighbor being infected.
}
\end{theorem}
\begin{IEEEproof}
See appendix.
\end{IEEEproof}

\begin{remark}[Lumped Node Optimal Infection Threshold]
{\rm
Note that in this work we show that such a threshold exist, and node``1'' will decide to connect or disconnect with the lumped node based on the chance of infection of the lumped node being below or above not this infection. However, we do not have a closed form solution for the optimal threshold $I^*$ for which~\eqref{eq:heuLP} is optimal, thus we rely MonteCarlo method to find the threshold or numerical optimization techniques to solve the problem. \oprocend
}
\end{remark}

\section{\textbf{Simulation}}

\subsection{Main Simulation} \label{se:simulation}
Consider a network of $5$ nodes $\mathcal{G}$, such that we are given the adjacency matrix $A$, the initial conditions are $\x^S(0)=[1,1,0.62,0.99,0.01]$, $\x^I(0)=[0,0,0.38,0.01,0.99]$, the rates of the model $\beta=0.2$, $\delta=1/10$, safety threshold of ($\threshold_1=0.7$), the parameter $c_{12}=10$, and the terminal time $T=11$ days. 

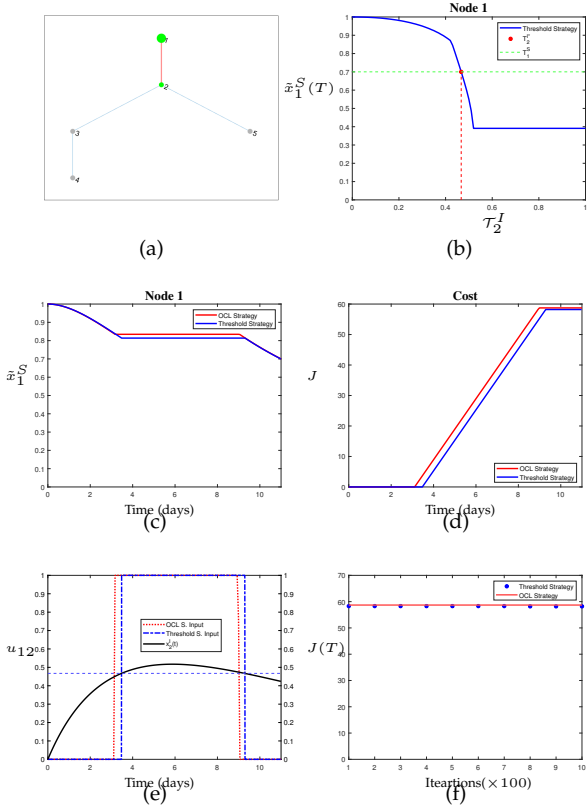
\begin{figure}[h] 
    \centering
\subfigure[]{\includegraphics[width=0.45\linewidth]{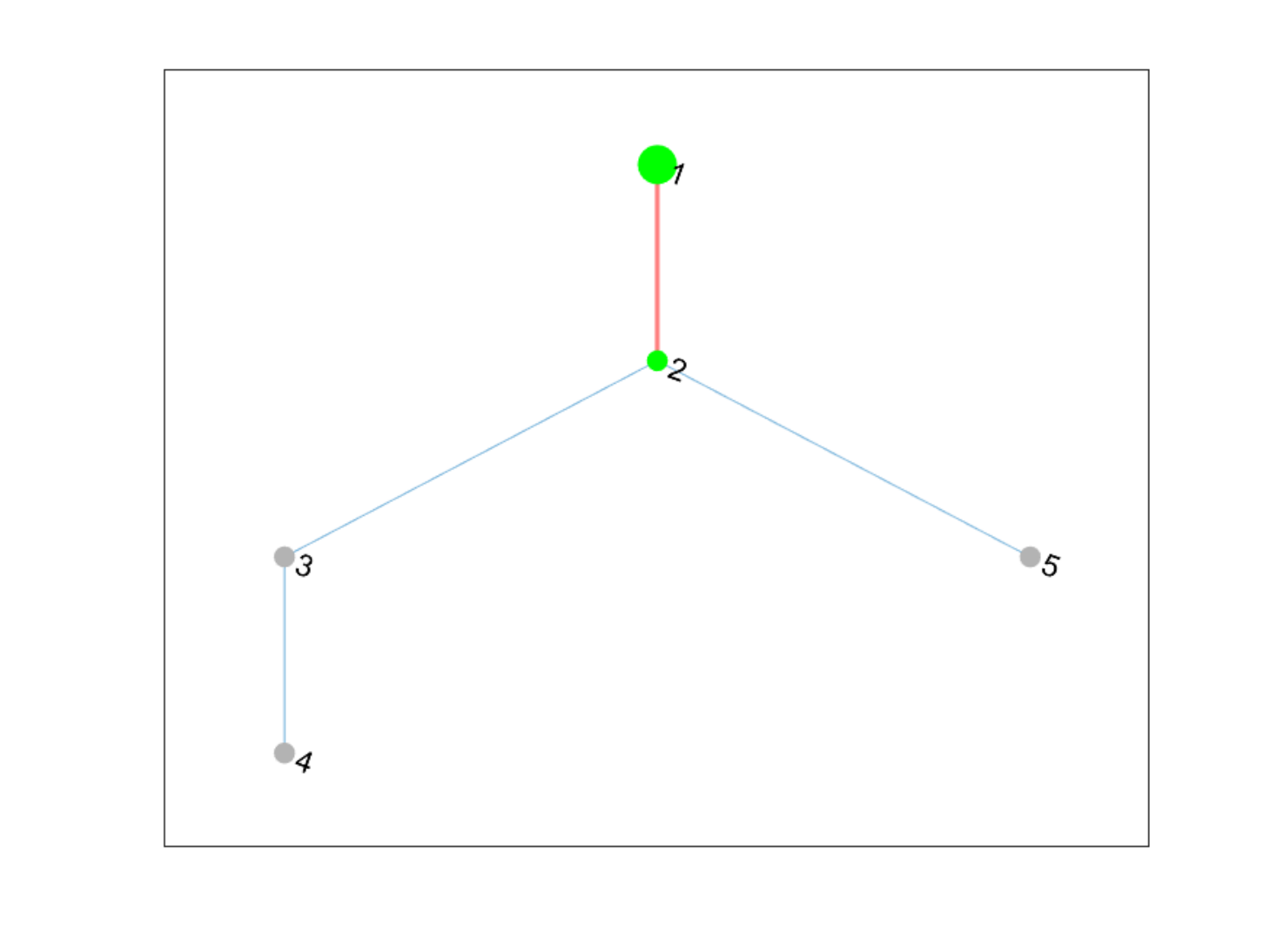}}
\subfigure[]{\includegraphics[width=0.45\linewidth]{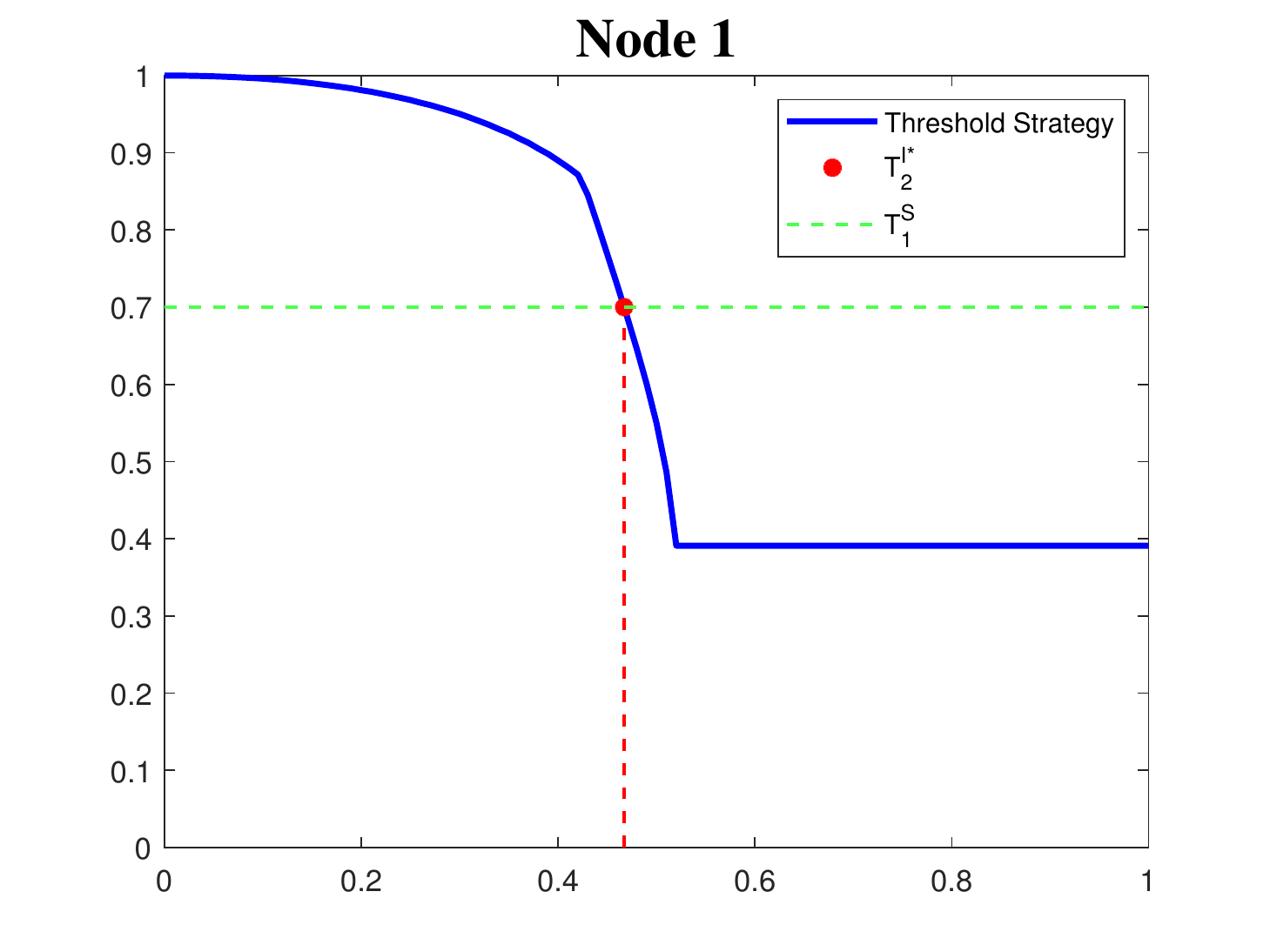}}
\tiny \put(-125,50){$\x_1^S(T)$}
\tiny \put(-50,-1){$\mathcal{T}_2^I$}\\
\subfigure[]{\includegraphics[width=0.45\linewidth]{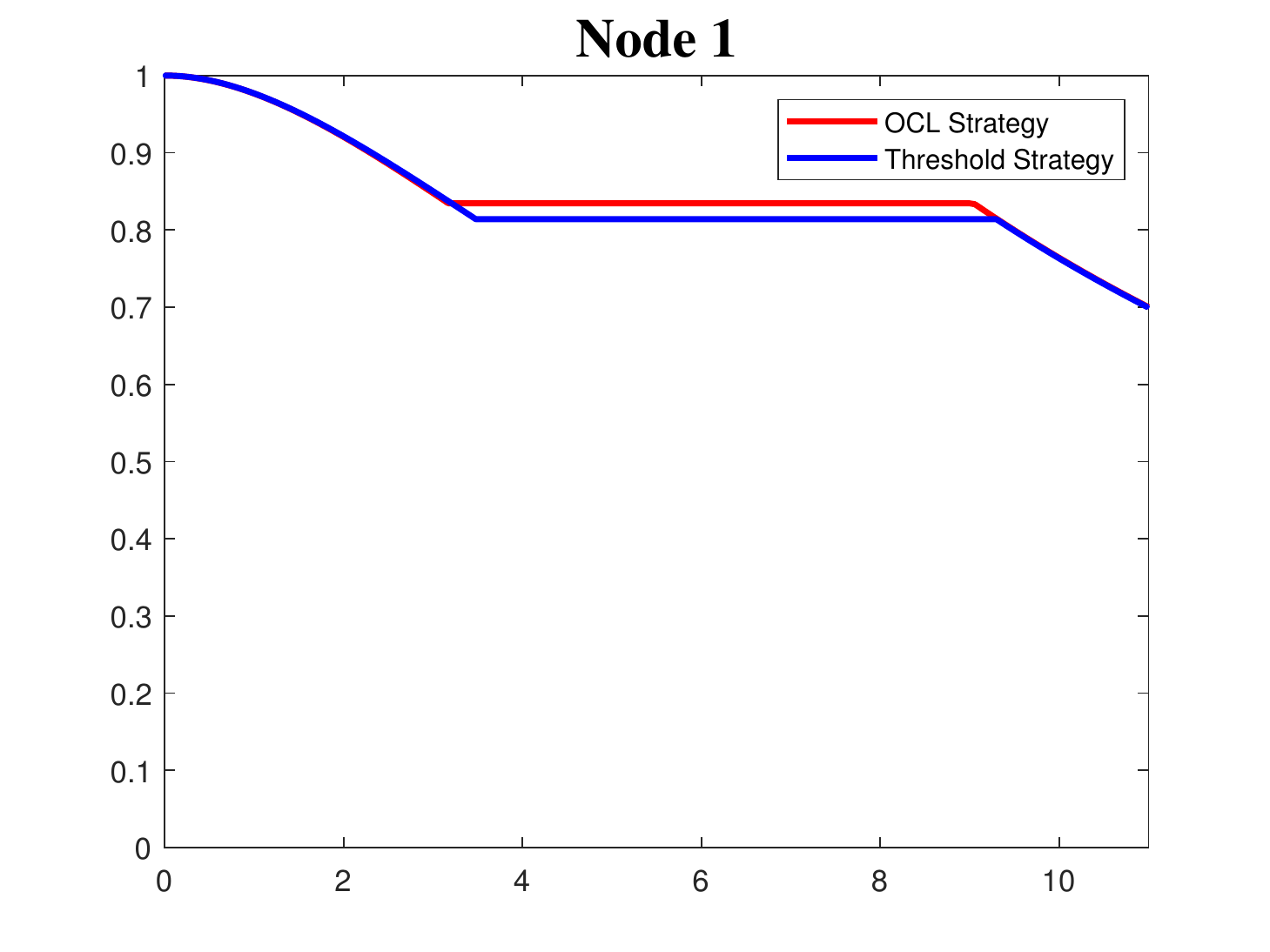}}
\tiny \put(-115,50){$\x_1^S$}
\tiny \put(-70,-1){Time (days)}
\subfigure[]{\includegraphics[width=0.45\linewidth]{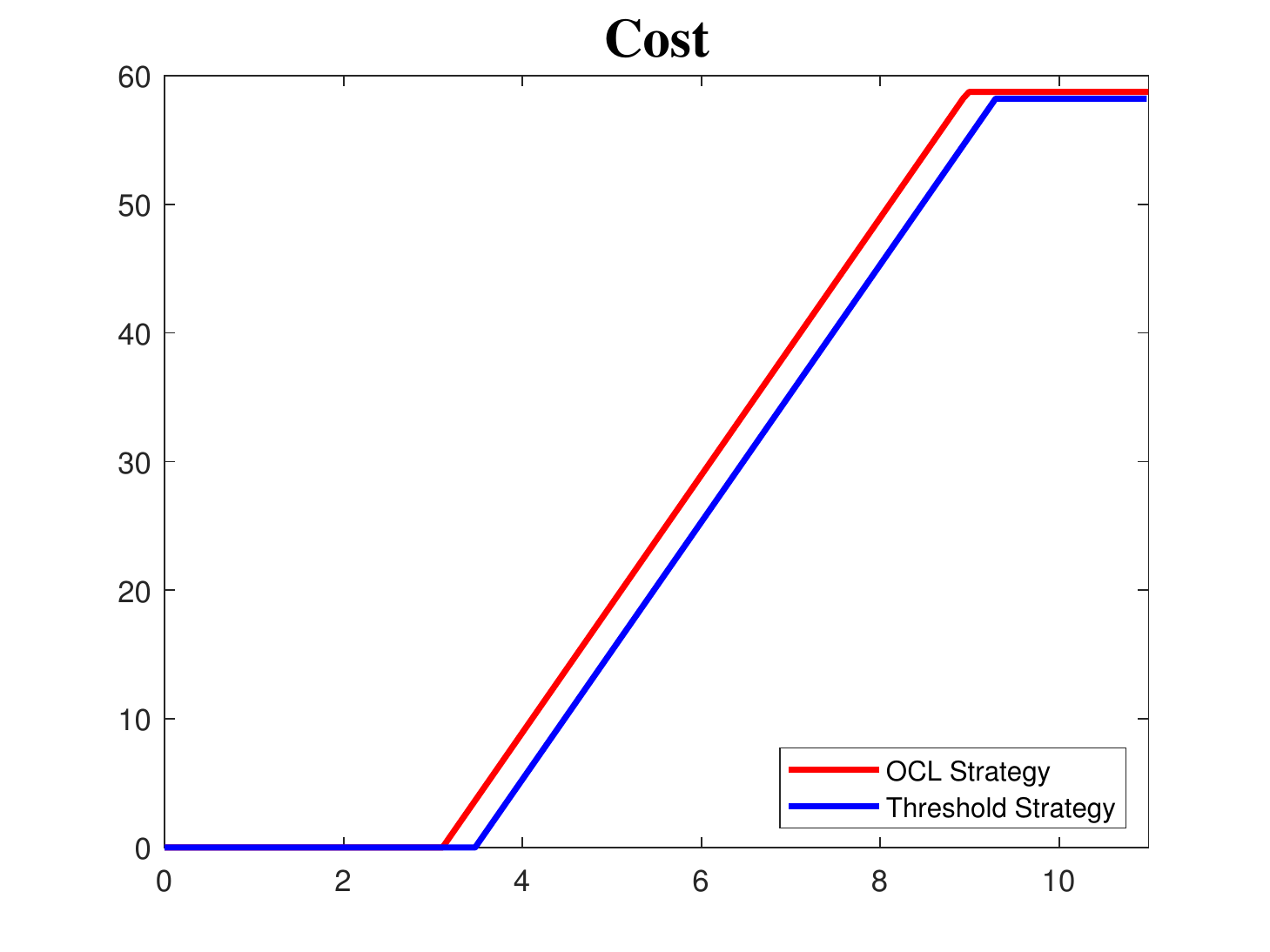}}
\tiny \put(-115,50){$J$}
\tiny \put(-70,-1){Time (days)} \\
\subfigure[]{\includegraphics[width=0.45\linewidth]{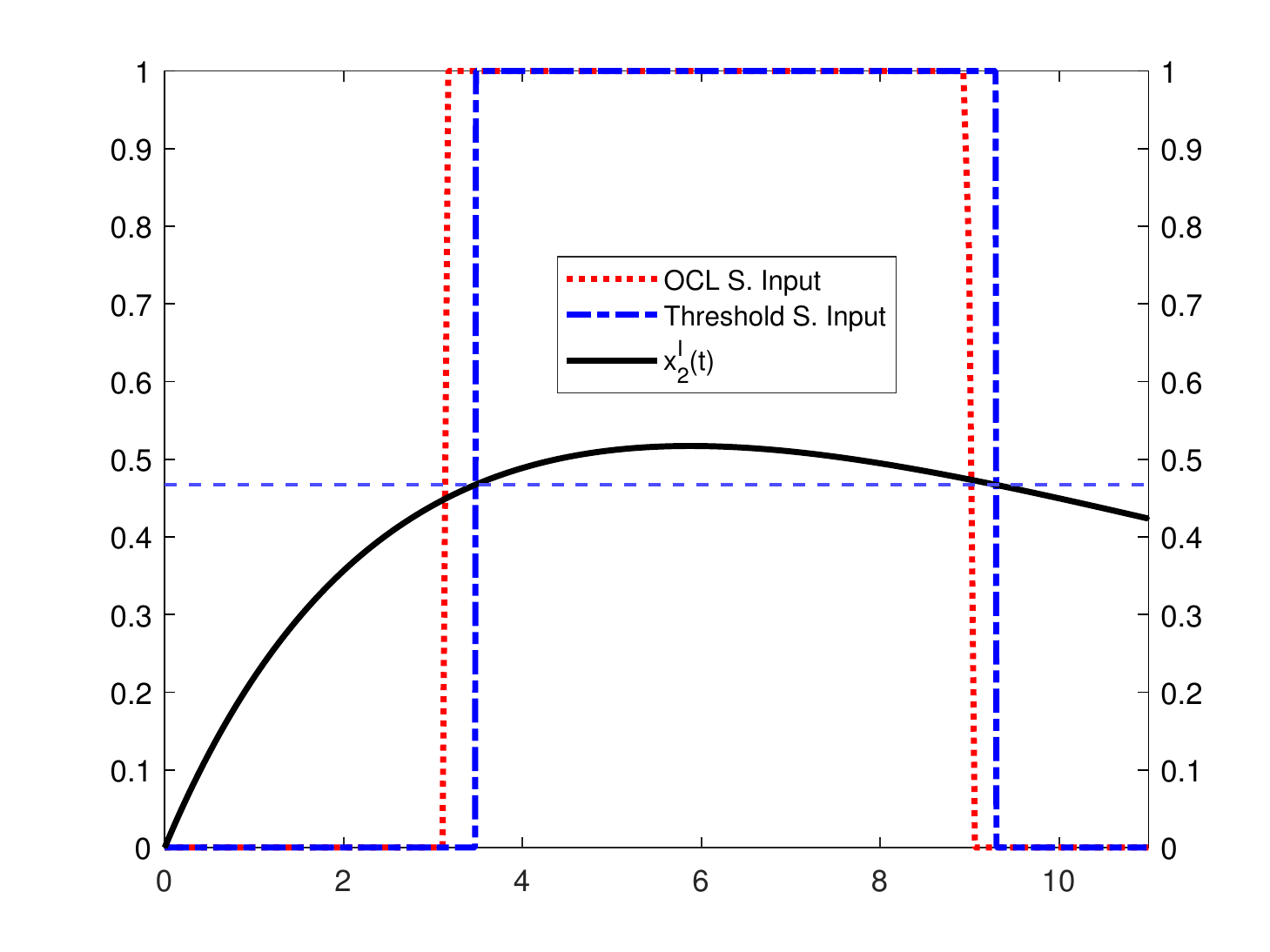}}
\tiny \put(-115,50){$u_{12}$}
\tiny \put(-70,-1){Time (days)}
\subfigure[]{\includegraphics[width=0.45\linewidth]{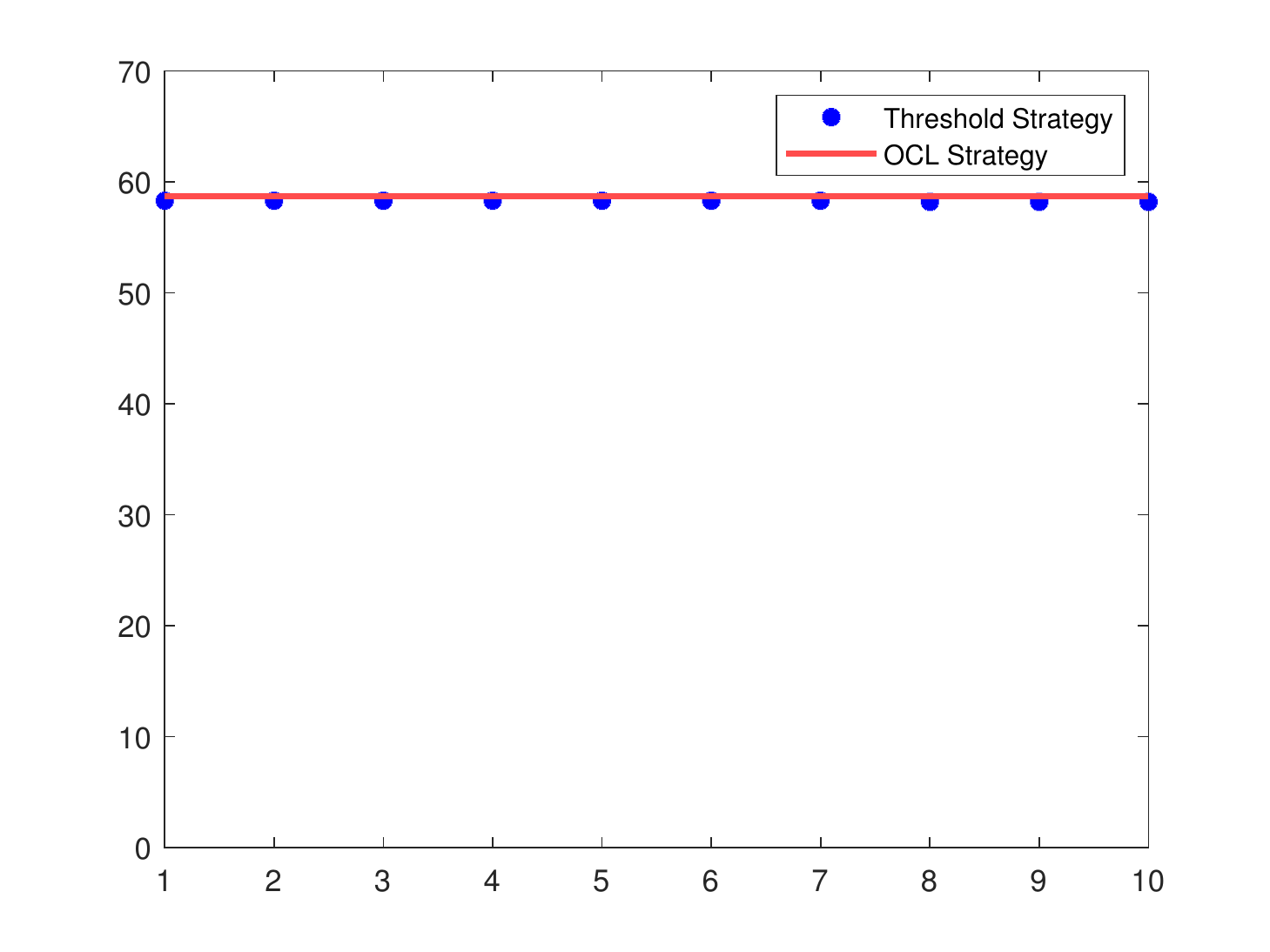}}
\tiny \put(-115,50){$J(T)$}
\tiny \put(-70,-1){Iteartions($\times 100$)}
    \caption{ (a) Network of 5 nodes, Node $1$ is enlarged and the red edge connects node $1$ to its direct neighbor. (b) Probability of node``1'' being susceptible at terminal time ($\x_1^S(T)$) against different values of $\mathcal{T}_2^I\in [0,1]$. The red dot is $\mathcal{T}_2^{I*}=0.4672$. (c) Probability of node``1'' being susceptible over time ($\x_1^S(t)$), for ''OCL`` and ``Threshold'' strategies. (d) the controller cost over time ($J(t)$). (e) the control, $u_{12}(t)$ (left) and the probability of node``2'' being infected (right) over time ($\x_2^I(t)$).}
    \label{fig:node1safe}
\end{figure}

In the first numerical simulation, the graph of $A$ can be visually represented as in Figure~\ref{fig:node1safe} (a). Note that the nodes colors in the graph represent the initial state of the nodes. Green color represents a Susceptible state, and grey color represent nodes that aren't $100\%$ susceptible or infected, i.e., refer the initial condition (distribution) of the nodes. Further, the enlarged node (i.e., node``1'') is the node that we keep its safety above a given safety threshold ($\threshold_1$). 

At the beginning, we find the optimal $\mathcal{T}^{I*}_2$, by using MonteCarlo method $\infectedthreshold_2 \in[0,1]$ as can be seen in~Figure \ref{fig:node1safe} (b). Such that the red dot is the optimal infection threshold of node ``2'', i.e., $\mathcal{T}^{I*}_2=0.4672$.


Since now we know $\mathcal{T}^{I*}_2$, therefore, we implement the ``Threshold Strategy'', which is the strategy of the solution in \eqref{eq:heu}. Further, we compare ``Threshold Strategy''(blue line) against ``OCL Strategy''(red line).  The ``OCL Strategy'' is the numerical optimization technique to solve Problem~\eqref{pr:CT-MFA}, where we use OpenOCL (and CasADi) software packages.

The comparison can be seen in Figure~\ref{fig:node1safe} (c) Probability of node``1'' being susceptible ($\x_1^S(t)$) against time(days), (d) the cost $J(t)$ against time(days), and (e) the control, $u_{12}(t)$ (left) and the probability of node``2'' being infected ($\x_2^I(t)$) (right) against time(days). We also show Figure~\ref{fig:node1safe} (f) the terminal cost vs the number of iteration when using MC method, where it can be observed that as the number of iterations increase the solution converges to the minimum cost.

\begin{figure}[h]
	\centering
	\includegraphics[width=0.9\linewidth]{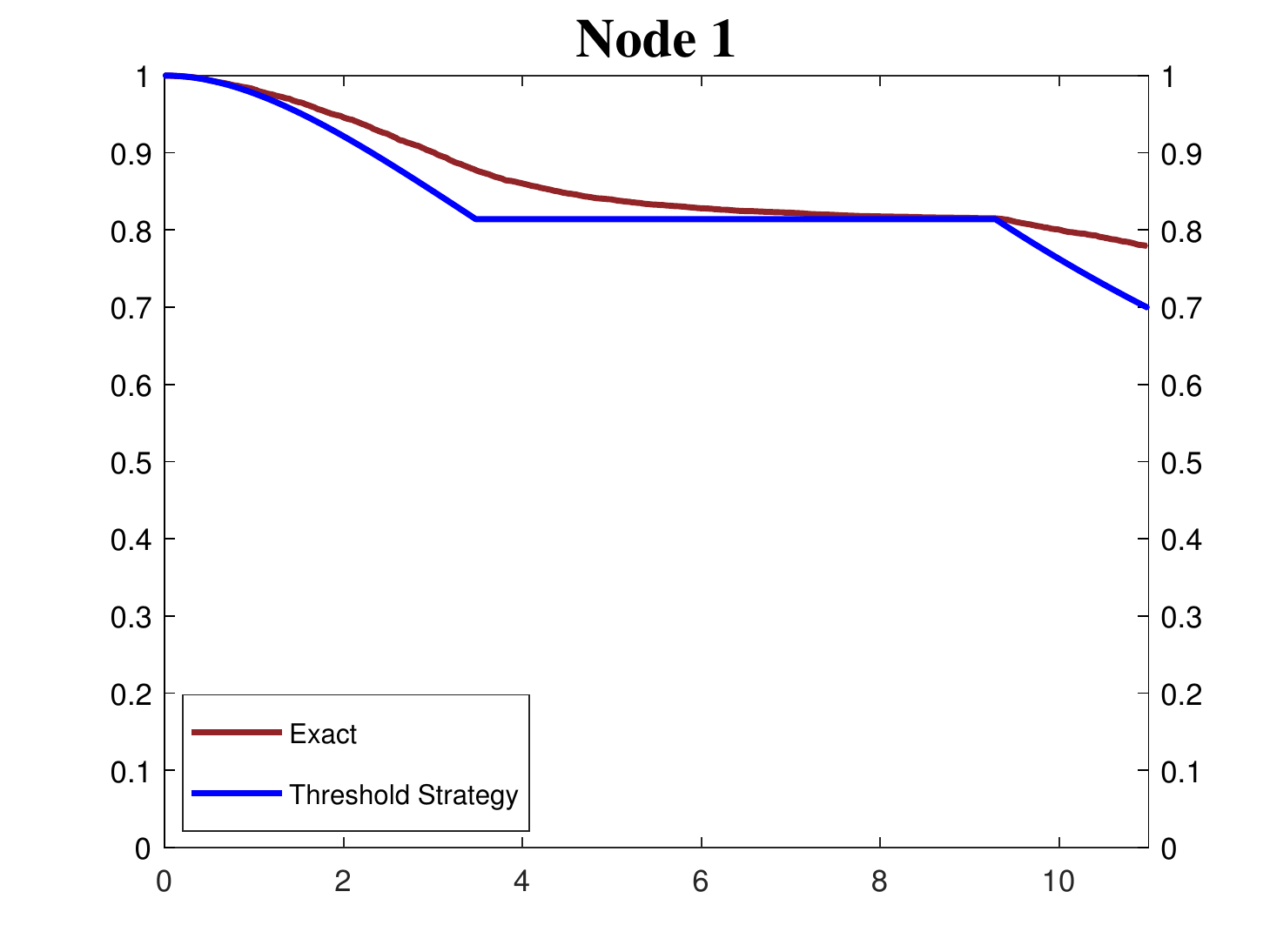}
	 \put(-220,80){$\x_1^S$}
	 \put(-3,80){$x_1^S$}
	 \put(-120,3){Time (days)}
	\caption{Showing the solution to the relaxed problem satisfy the safety threshold on the original stochastic problem. Brown line (right) is the exact probability of node ``1'' being susceptible ($x_1^S(t)$) over time when substituting the solution from the relaxed problem ``Threshold Strategy''. Blue line (left) is the probability of node ``1'' being susceptible when using ``Threshold Strategy''. }
\label{fig:exactVsThreshold}
\end{figure}

We would like also to show the results mentioned in Theorem~\ref{th:Inequality}, 
which show the effectiveness of solving problem~\ref{pr:CT-MFA} on the original problem~\ref{pr:SC}. Figure~\ref{fig:exactVsThreshold} verifies the results in Theorem~\eqref{th:Inequality}. We clearly can see the bound on the $x_1^S(t)$ by $\x_1^S(t)$, and further we can see $x_1^S(T)>\x_1^S(T)$.

\begin{figure}[h] 
    \centering
    \subfigure[]{\includegraphics[width=0.45\linewidth]{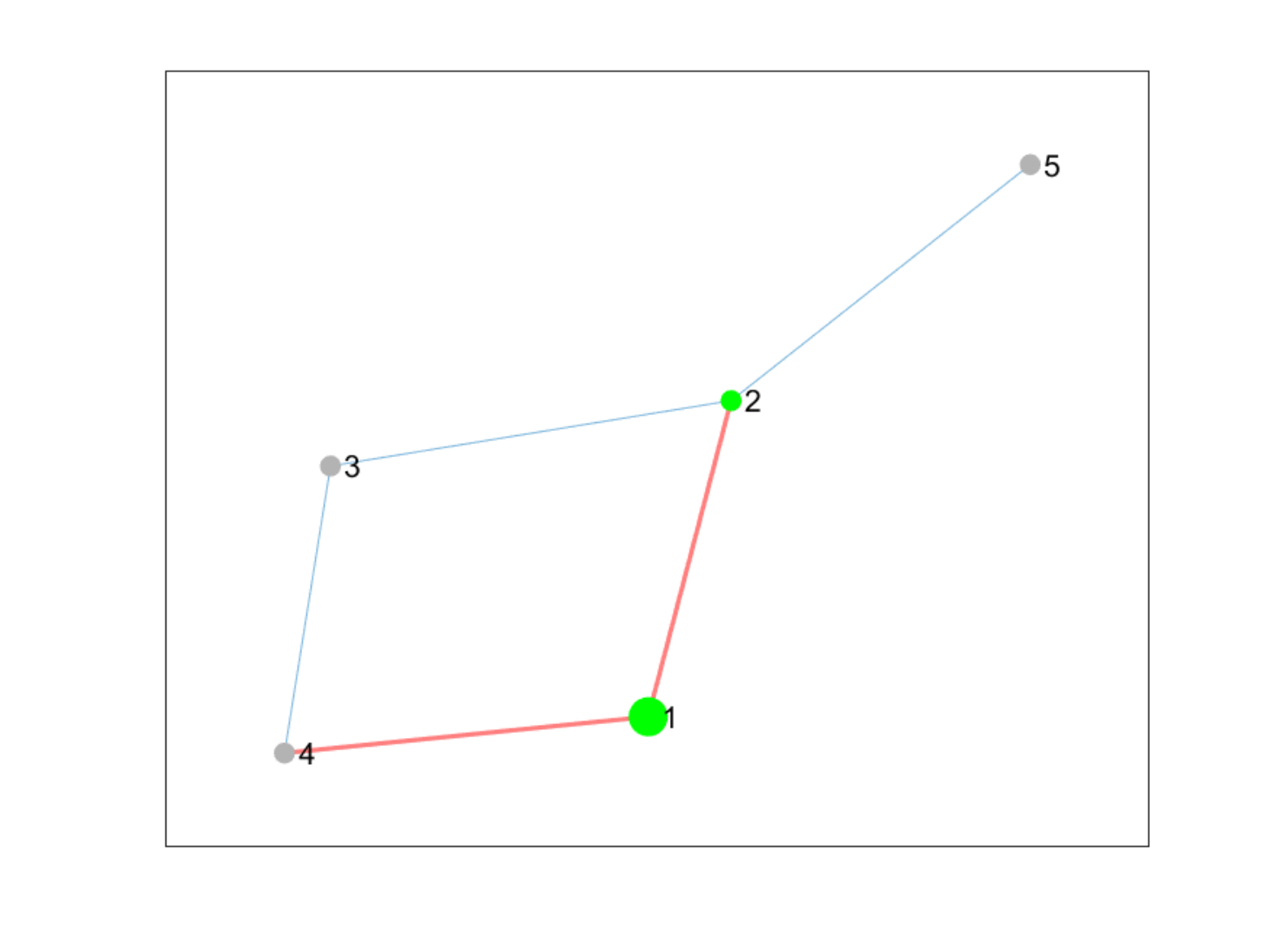}}
\tiny \put(-40,30){$c_{12}=10$}
\tiny \put(-80,10){$c_{14}=1$}
\subfigure[]{\includegraphics[width=0.45\linewidth]{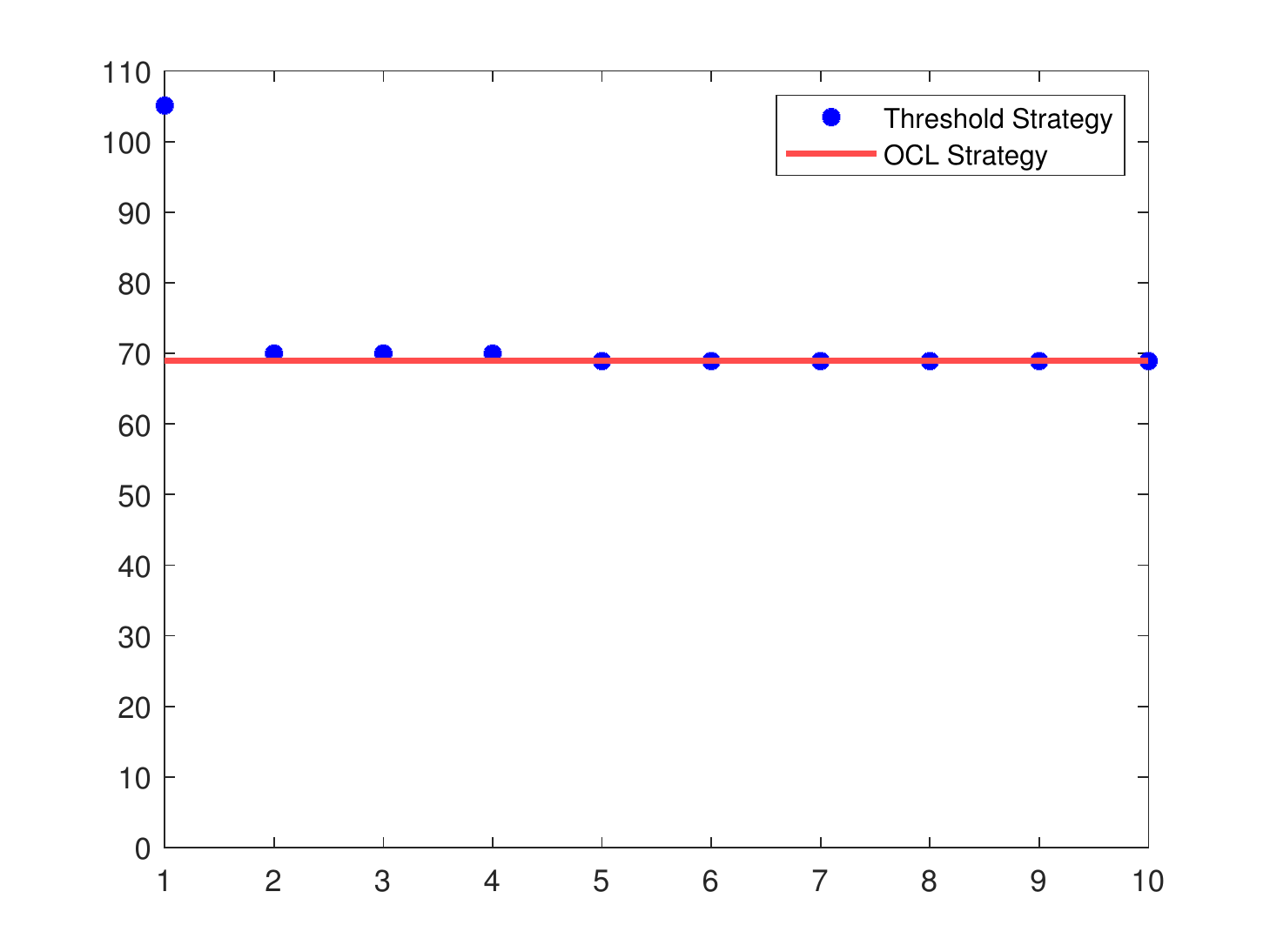}}
\tiny \put(-120,50){$J(T)$}
\tiny \put(-70,-1){Iteartions($\times 1000$)}\\
\subfigure[]{\includegraphics[width=0.45\linewidth]{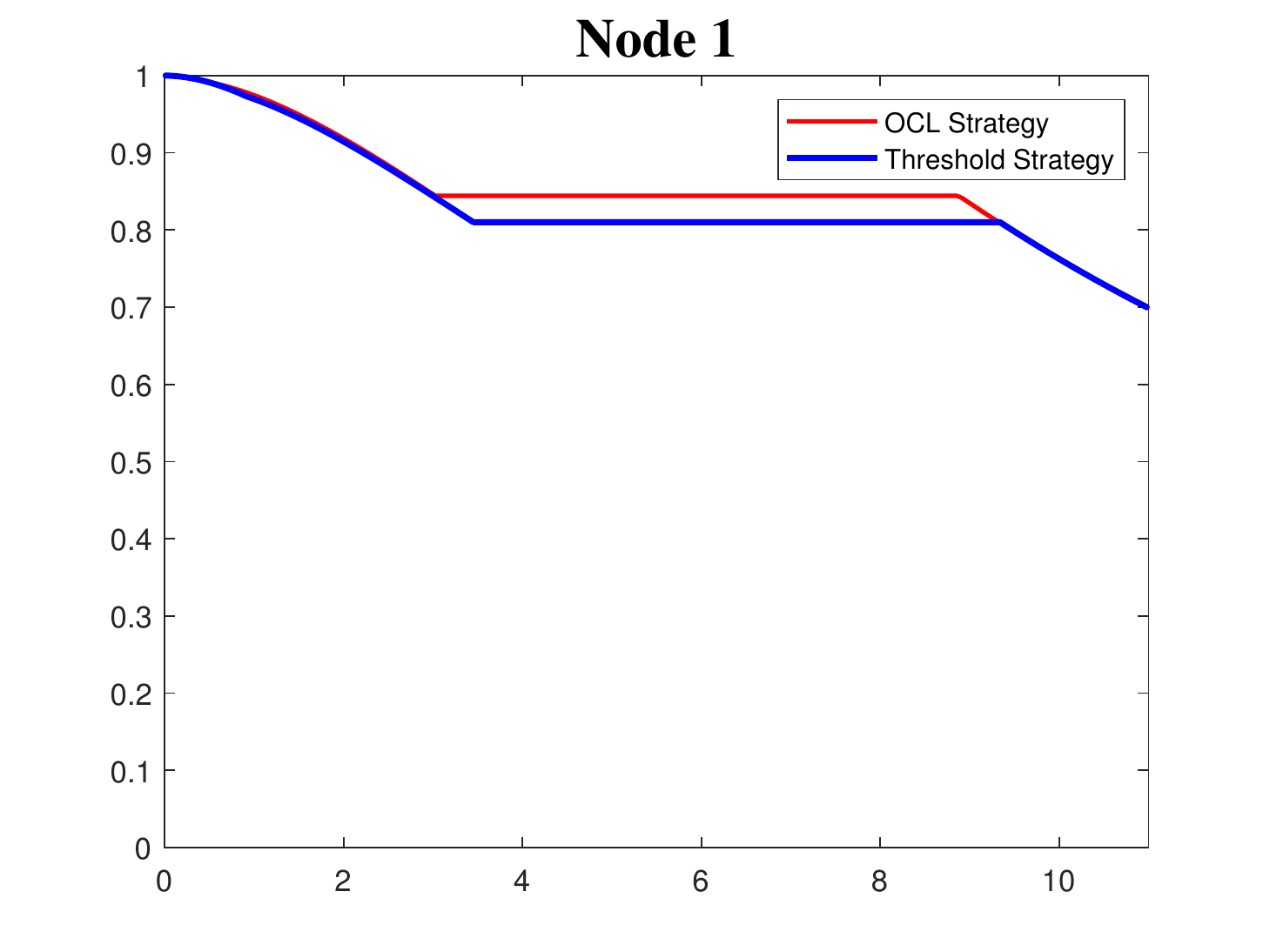}}
\tiny \put(-115,50){$\x_1^S$}
\tiny \put(-70,1){Time (days)}
\subfigure[]{\includegraphics[width=0.45\linewidth]{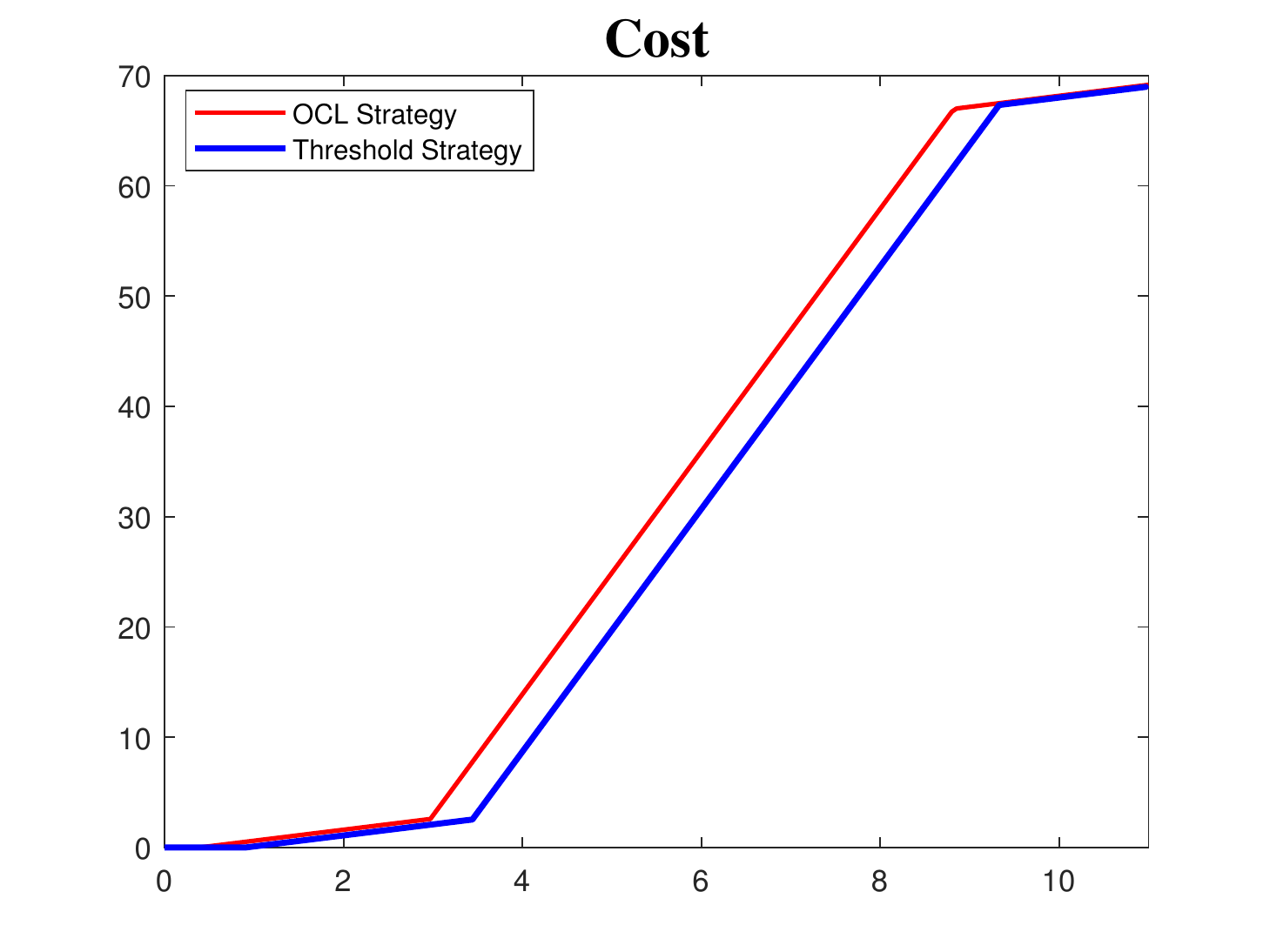}}
\tiny \put(-115,50){$J$}
\tiny \put(-70,1){Time (days)}\\
\subfigure[]{\includegraphics[width=0.45\linewidth]{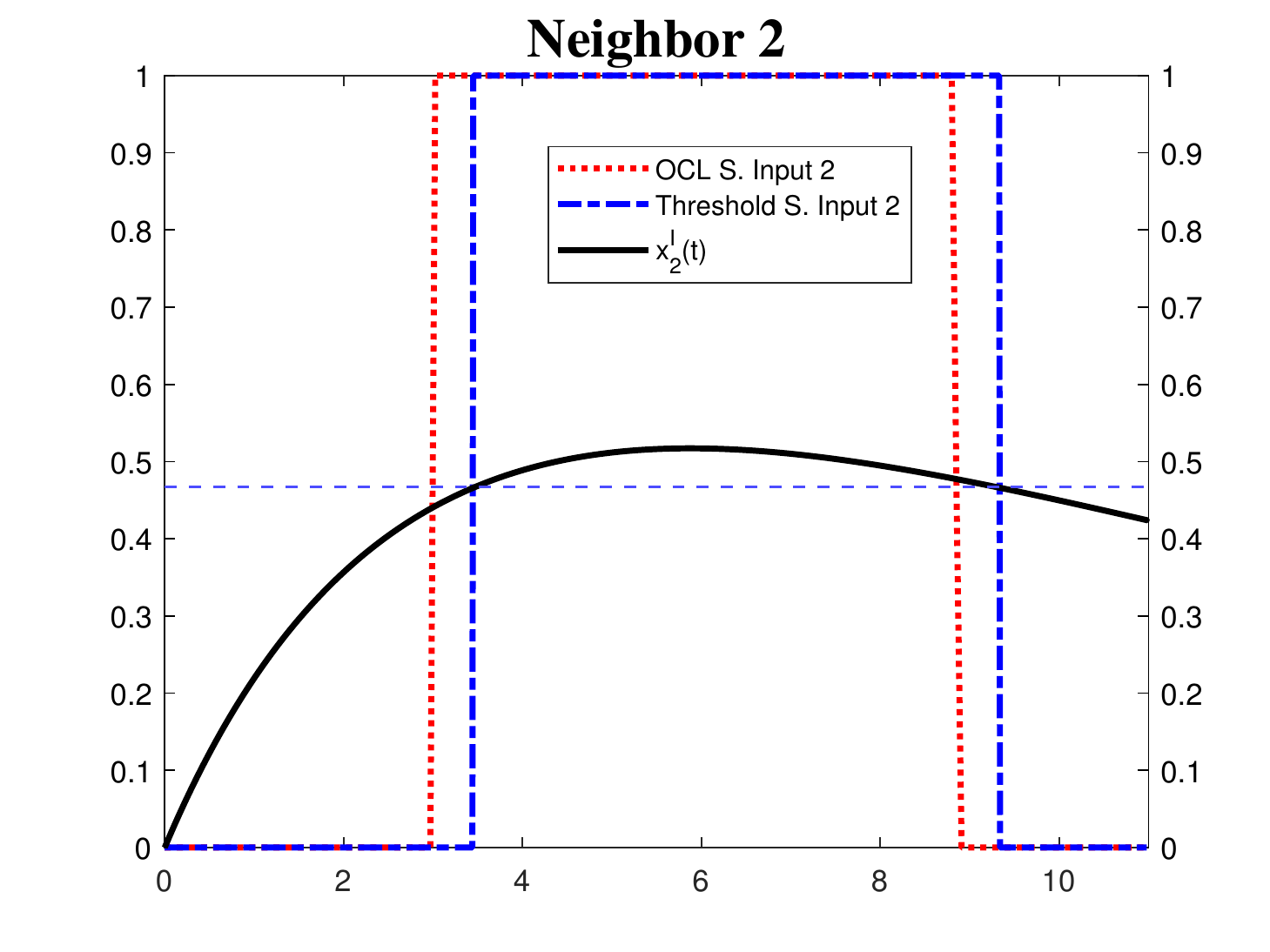}}
\tiny \put(-115,50){$u_{12}$}
\tiny \put(-50,1){Time (days)}
\subfigure[]{\includegraphics[width=0.45\linewidth]{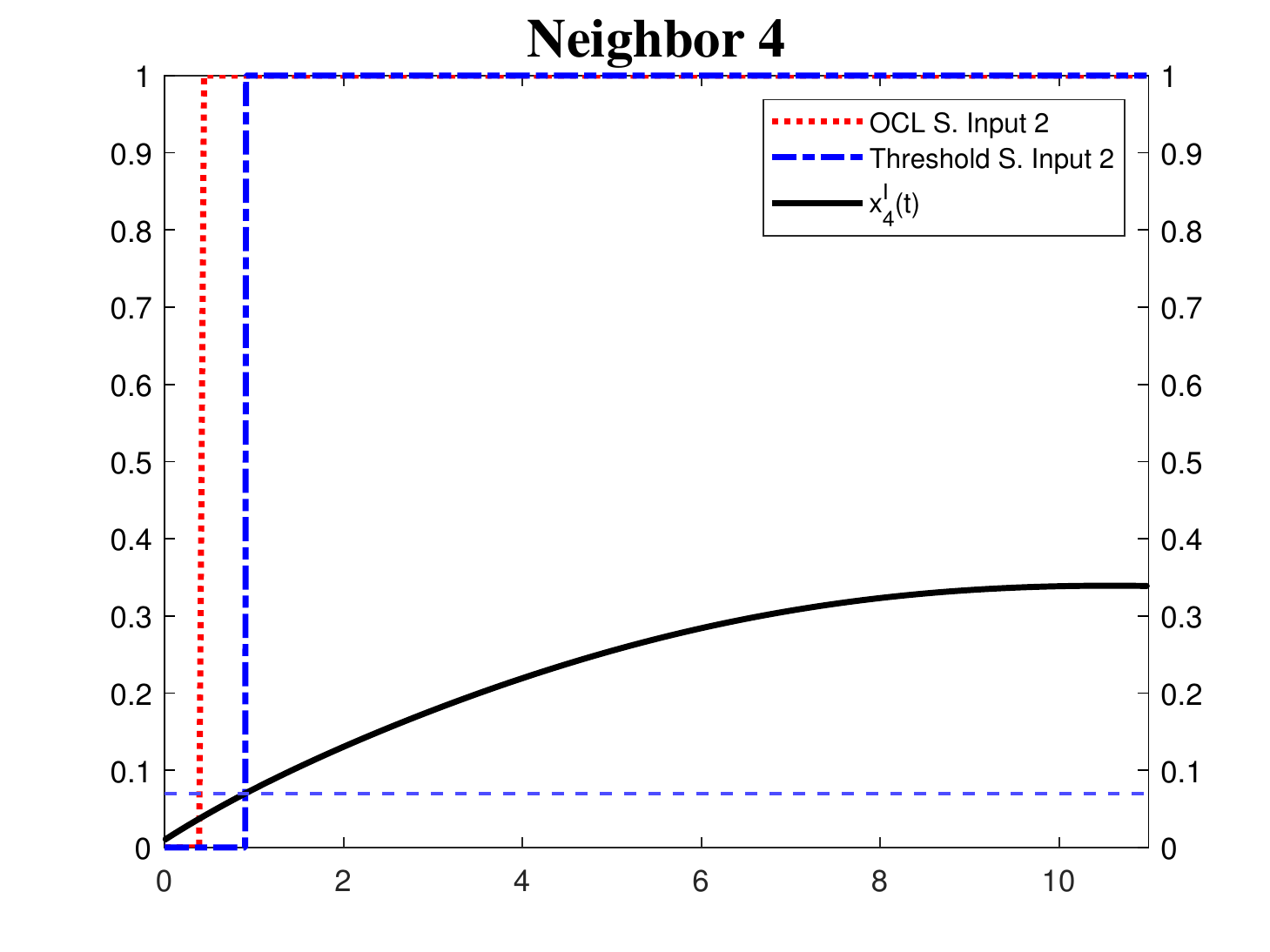}}
\tiny \put(-115,50){$u_{14}$}
\tiny \put(-70,1){Time (days)}
\caption{(a) Network of 5 nodes. Node $1$ is enlarged, which has two direct neighbors. (b) Number of iterations needed to give an accurate terminal cost when using MC method (c) Probability of node``1'' being susceptible over time ($\x_1^S(t)$), for ''OCL`` and ``Threshold'' strategies. (d) The controller cost over time ($J(t)$). (e) The control input $u_{12}(t)$ (left), and probability of node ``2'' being infected(left) over time ($\x^I_2(t)$). (f) The control input $u_{14}(t)$ (left), and probability of node``4'' being infected (right) over time $\x^I_4(t)$.} 
    \label{fig:5nodes}
\end{figure}

Secondly, as a second numerical simulation, we would like to show a case when node ``1'' has more than one neighbor. We use the same given initial conditions, and the new $A$ matrix is given and the graph is shown in~ Figure~\eqref{fig:5nodes} (a). The $c_{12}=10$, and $c_{14}=1$ are given, which also can be seen in the figure. This means that disconnecting from node ``2'' will cost more than disconnecting from node ``4''. 
Using MonteCarlo (MC) method, we find $\mathcal{T}_2^{I*}=0.466$ and $\mathcal{T}_4^{I*}=0.07$. In Figure~\eqref{fig:5nodes} (b) we show the terminal cost vs the number of iteration when using MC method, where it can be observed that as the number of iterations increase the solution converges to the minimum cost. Then, we show the results of implementing ``Threshold Strategy'' vs ``OCL Strategy'' in Figure~\eqref{fig:5nodes} (c) $\x_1^S(t)$ against time (days), (d) $J(t)$ against time(days), (e) $u_{12}(t)$ and $\x^I_2(t)$ against time, (f) $u_{14}(t)$ and $\x^I_4(t)$ against time. 


\begin{figure}[h] 
\centering
	\subfigure[]{\includegraphics[width=0.45\linewidth]{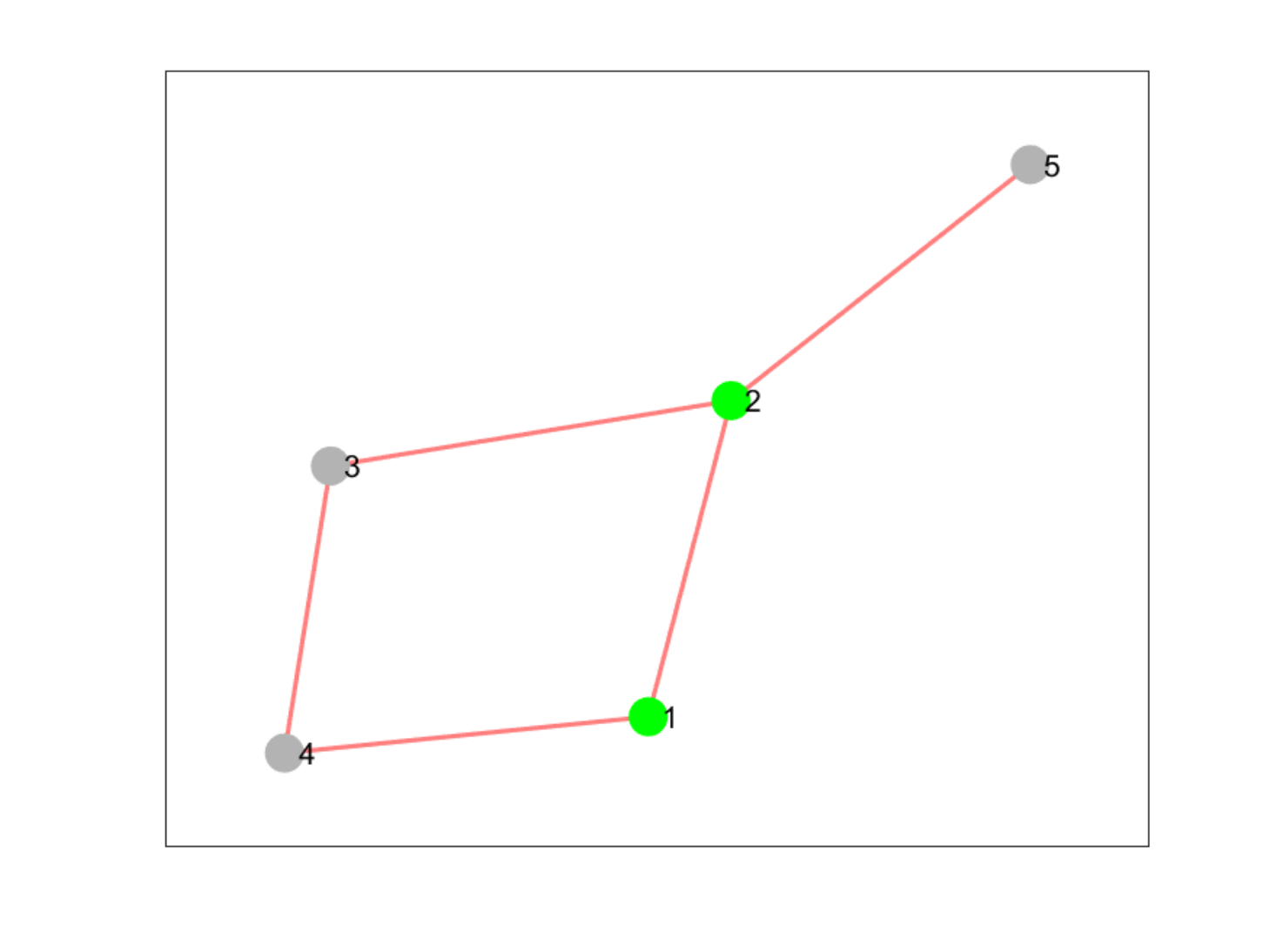}}\\
	\subfigure[]{\includegraphics[width=0.45\linewidth]{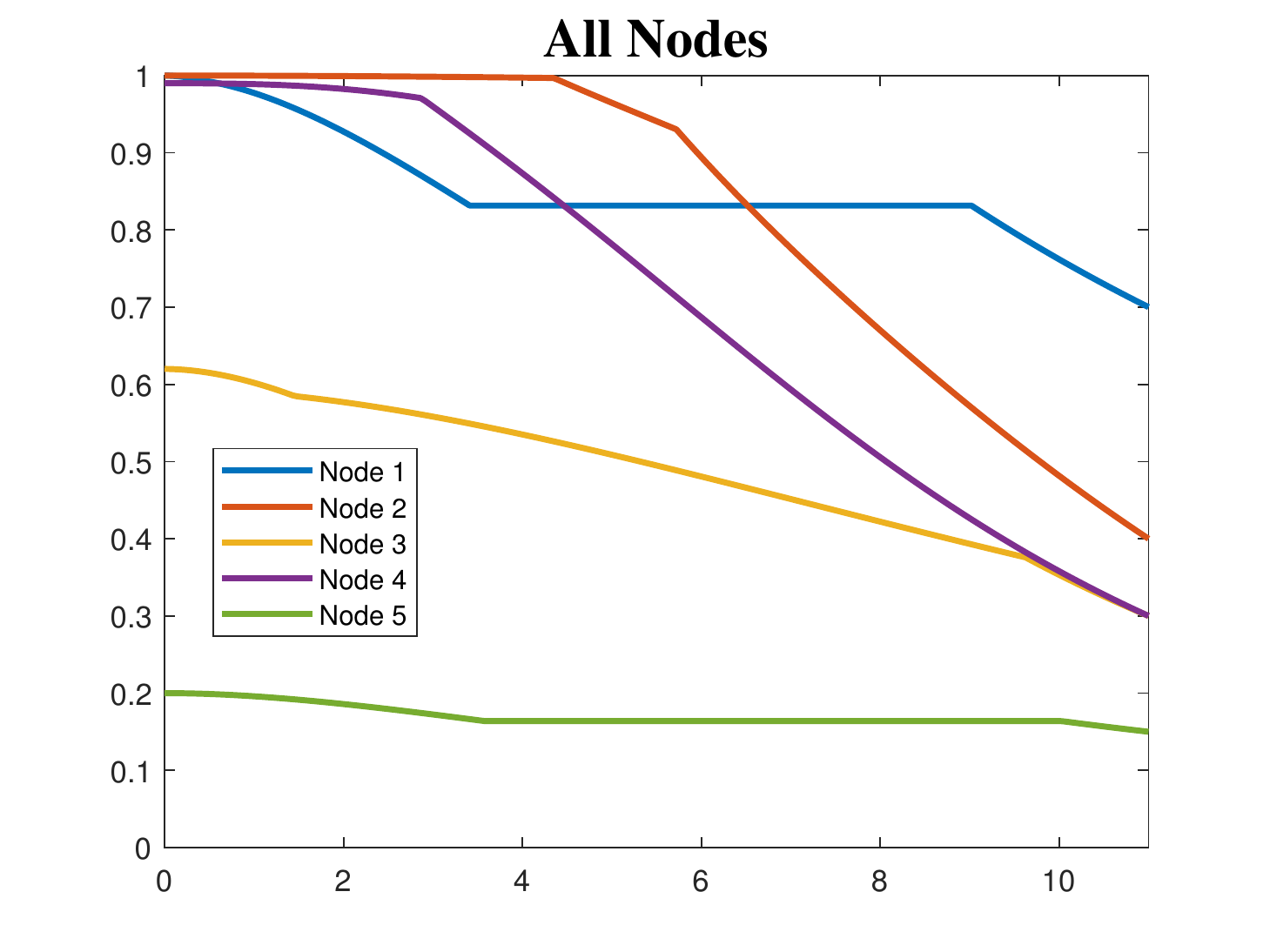}}
		\tiny \put(-115,50){$\x_i^S$}
		\tiny \put(-70,1){Time (days)}
	\subfigure[]{\includegraphics[width=0.45\linewidth]{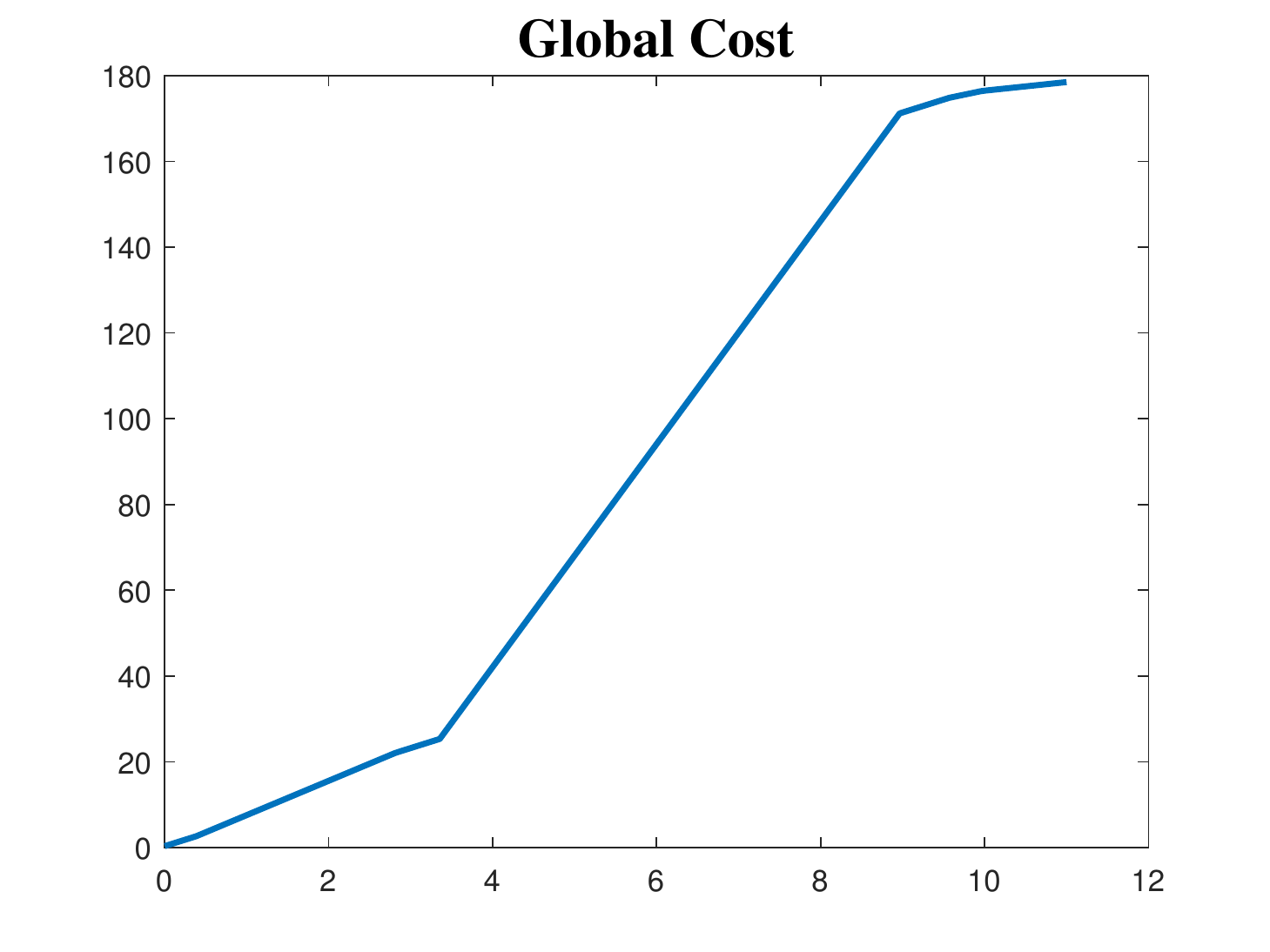}}
		\tiny \put(-115,50){$J_G$}
		\tiny \put(-70,1){Time (days)}
\caption{(a) Graph of 5 Node (b) Probability of all nodes being susceptible over time ($\x_i^S(t)$), where $i \in \V$, using 'OCL Strategy'. (c) the global cost of the controller over time $J_G(t)$.}
    \label{fig:globalcontrol}
\end{figure}

Thirdly, as a third numerical example, we show a case when all agents in the network is trying to keep themselves above some personal threshold, such that $\threshold=[0.7, 0.4, 0.3, 0.3, 0.15]$, and $c_{12}=c_{21}=10$, $c_{14}=c_{41}=1$, $c_{23}=c_{32}=1$, $c_{34}=c_{43}=1$, and $c_{25}=c_{52}=1$. Such that in this example every individual is running his local strategy based on the given information. The question would someone ask is what would the global cost for the entire network in this case be. Figure~\ref{fig:globalcontrol}, shows the probability of all nodes being susceptible, and the global cost over time. As we can see from Figure~\ref{fig:globalcontrol} (b), each need has it's own personal threshold, that is different from others in the network.

\begin{figure}[h] 
\centering
	\subfigure[]{\includegraphics[width=0.5\linewidth]{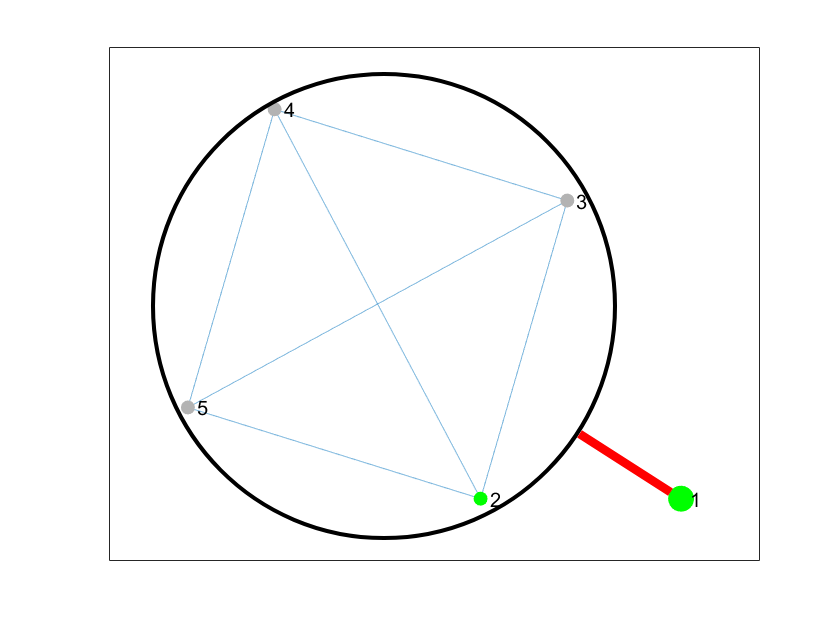}}
		\tiny \put(-35,30){$c_{av}=22$}\\
	\subfigure[]{\includegraphics[width=0.45\linewidth]{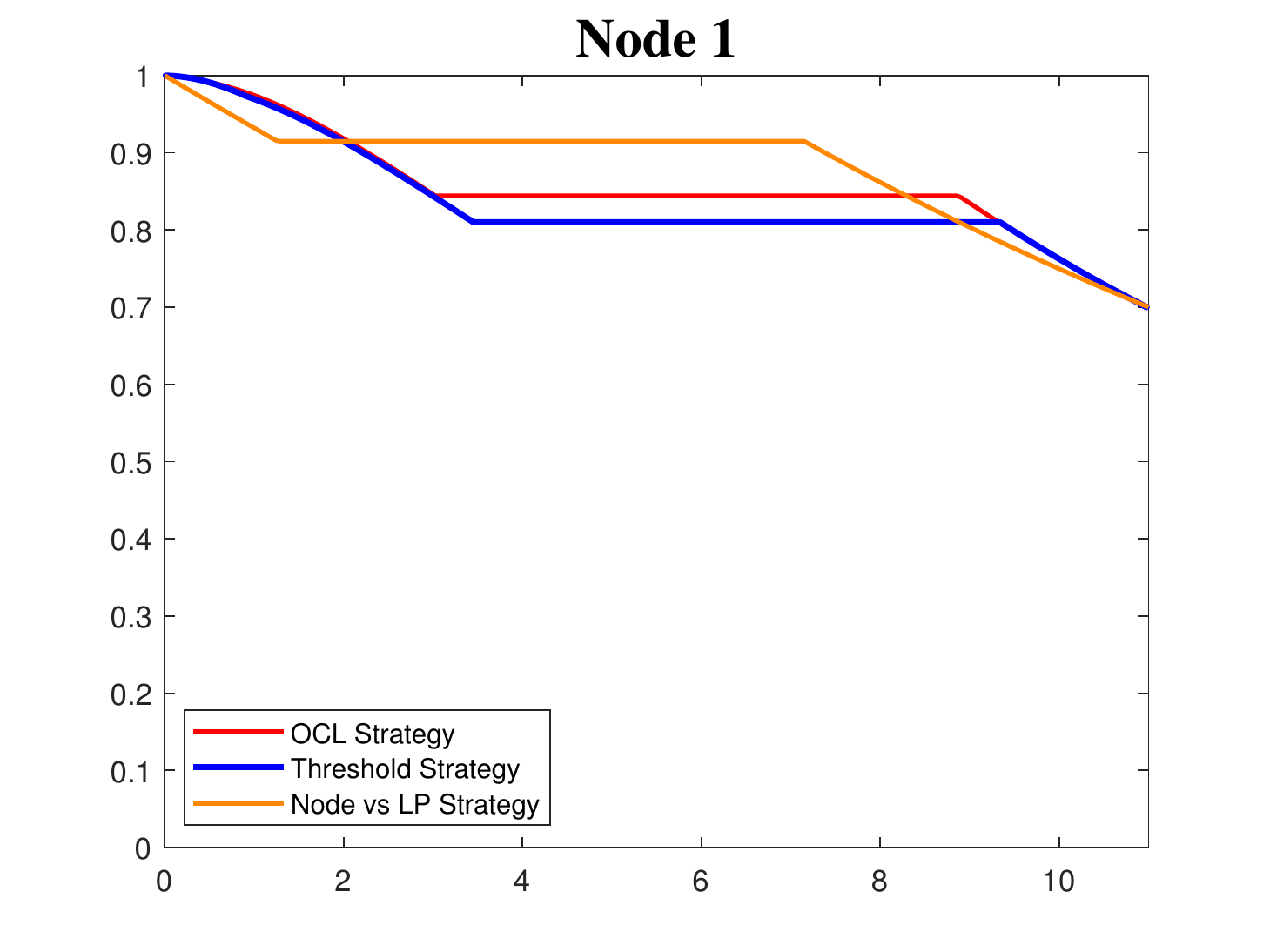}}
		\tiny \put(-115,50){$\x_1^S$}
		\tiny \put(-70,1){Time (days)}
	\subfigure[]{\includegraphics[width=0.45\linewidth]{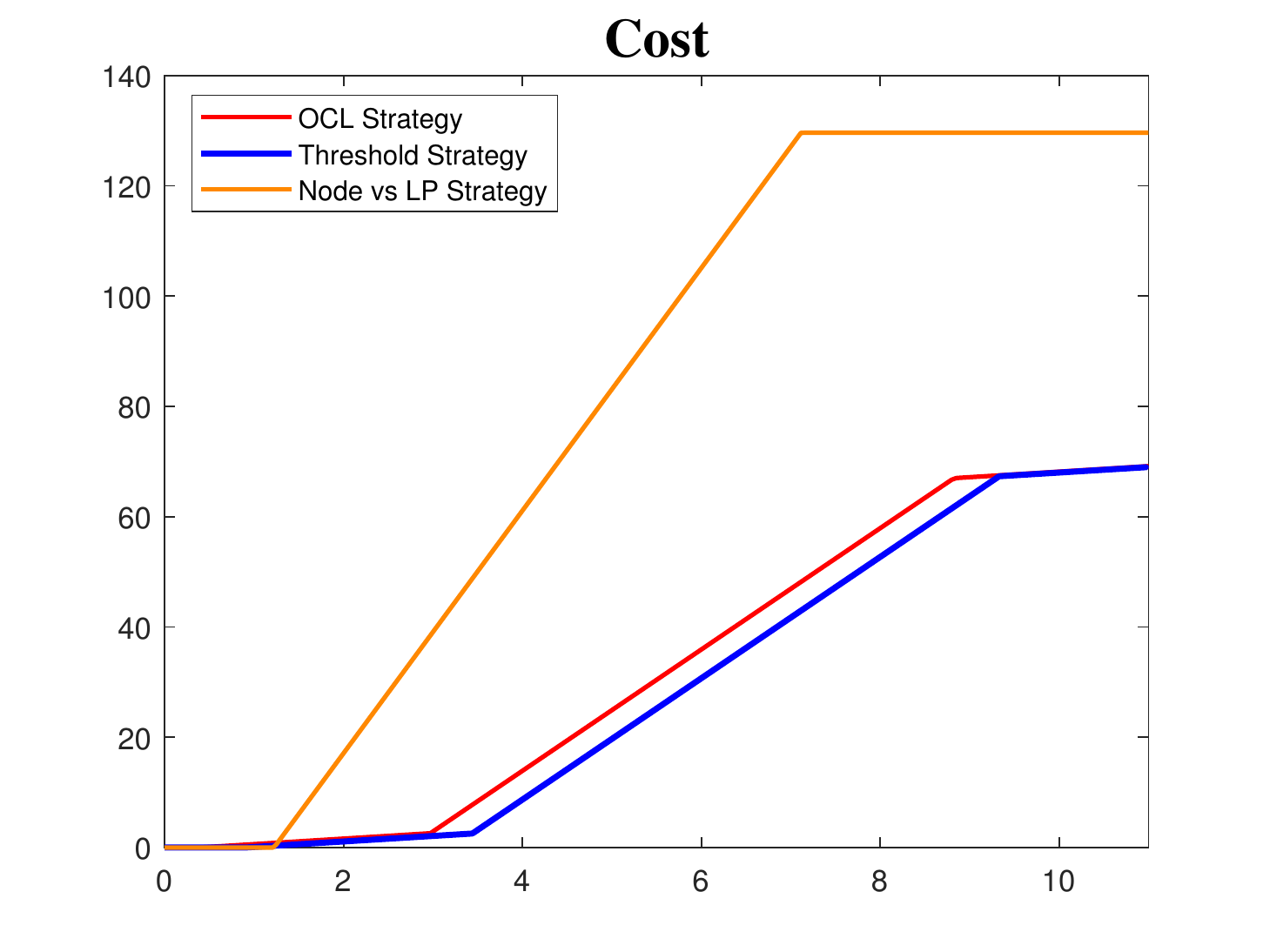}}
		\tiny \put(-110,50){$J$}
		\tiny \put(-70,1){Time (days)}
\caption{(a) graph of Node``1'' interacting with Lumped Node, inside the lumped node a complete graph is considered (b) Probability of node``1'' being susceptible over time ($\x_1^S(t)$), for different strategies. (c) the total cost of the controller over time $J(t)$.}
    \label{fig:n-vs-LP}
\end{figure}

Lastly, we compare the ``Threshold Strategy'' and ``OCL Strategy'' with ``Node vs LP Strategy'' (orange line). This last strategy suggest the well mixed homogeneous network. Such that node``1'' interactions can be lumped into one link,i.e., one control input, and other nodes on the graph can be lumped into one node by assuming a complete graph. The lumped graph can be seen in Figure~\eqref{fig:n-vs-LP} (a), where $S_{LP}=\frac{1}{4}\sum_{i=2}^5 \x_i^S(t)=0.6550$,  $I_{LP}=\frac{1}{4}\sum_{i=2}^5 \x_i^I(t)=0.345$, and $c_{av}=4\frac{c_{12}+c_{14}}{2}$.
The results of the comparison of the three strategies illustrated in Figure~\ref{fig:n-vs-LP} (b) and (c).


\subsection{Time Complexity}
\begin{figure*}[h] 
\centering
	\subfigure[]{\includegraphics[width=0.3\linewidth]{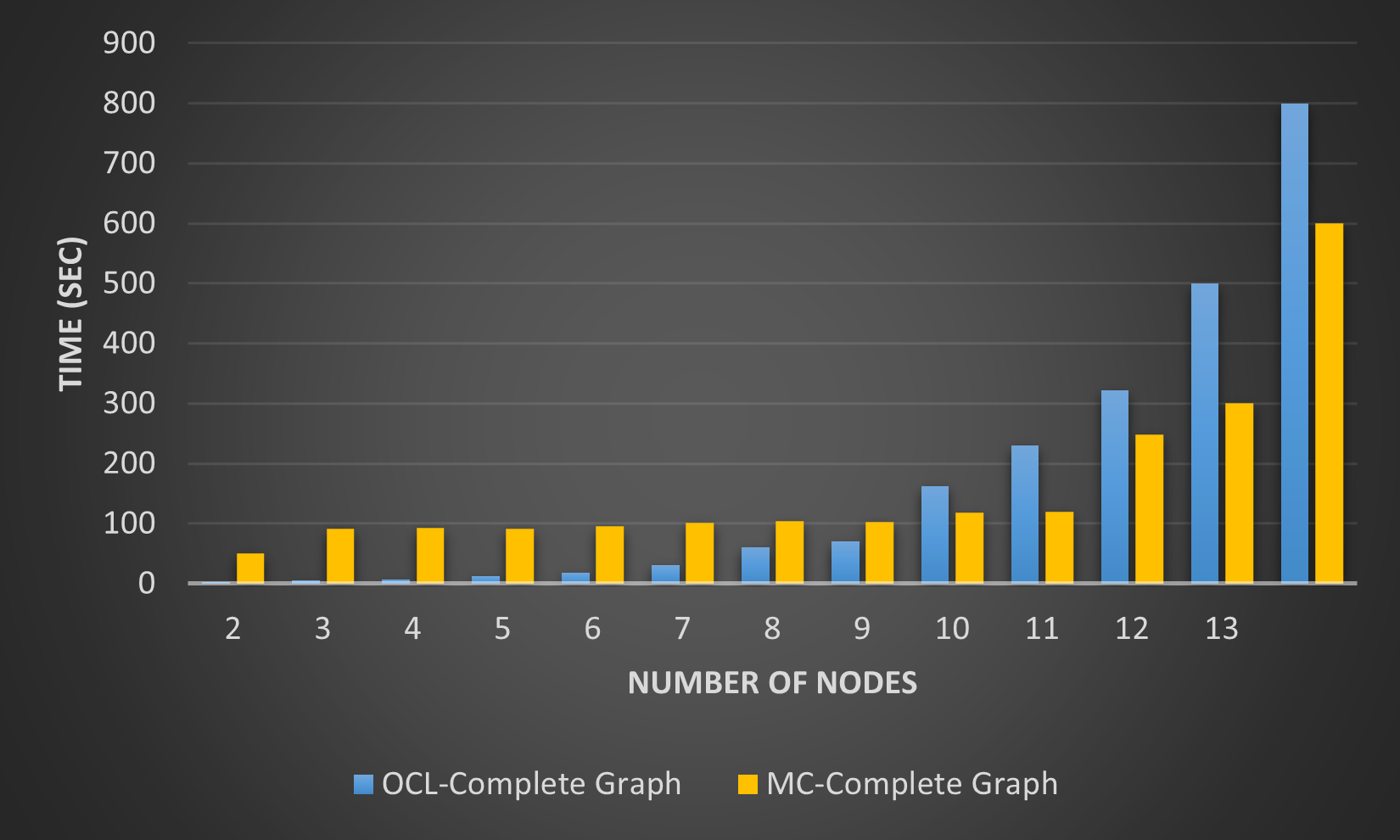}}
	\subfigure[]{\includegraphics[width=0.3\linewidth]{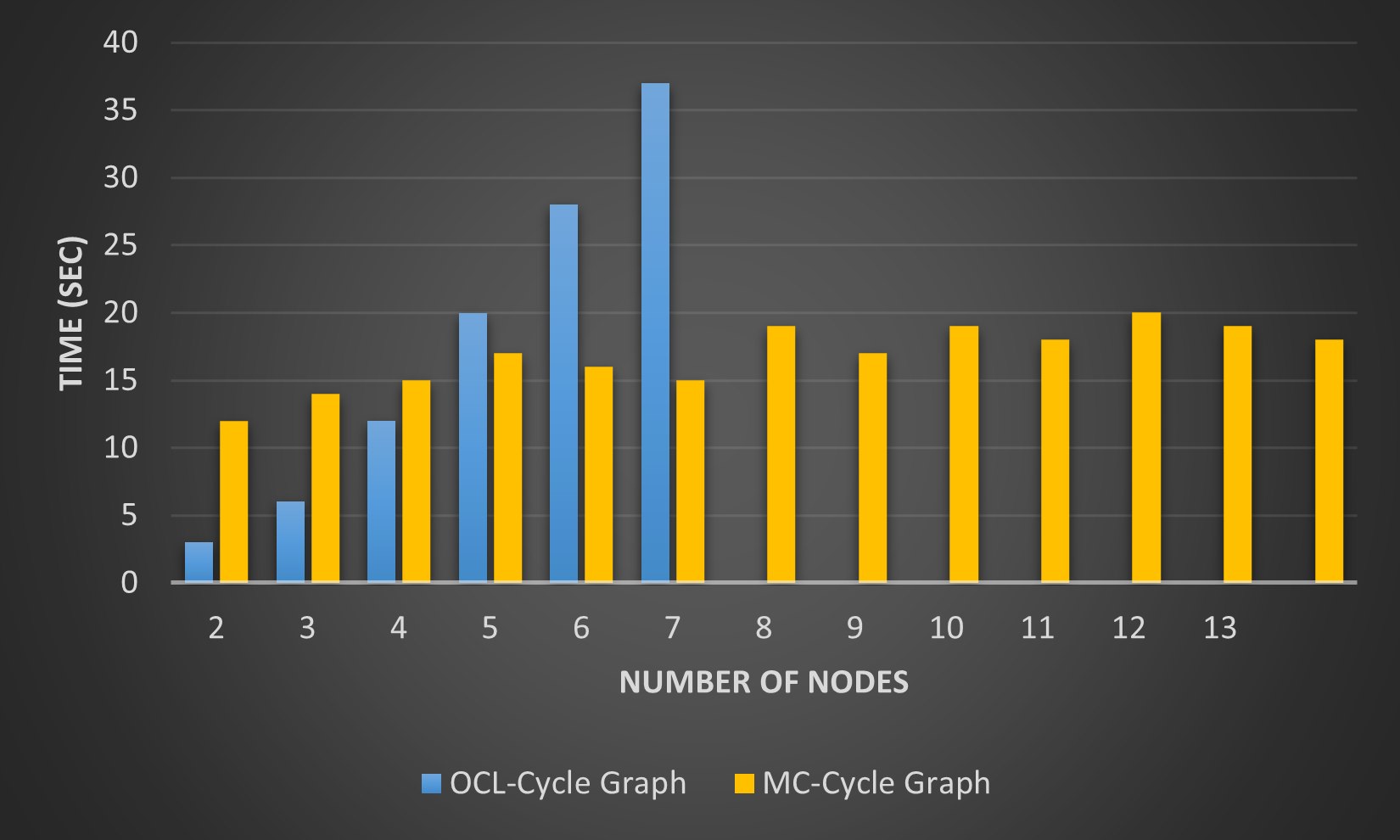}}
	\subfigure[]{\includegraphics[width=0.3\linewidth]{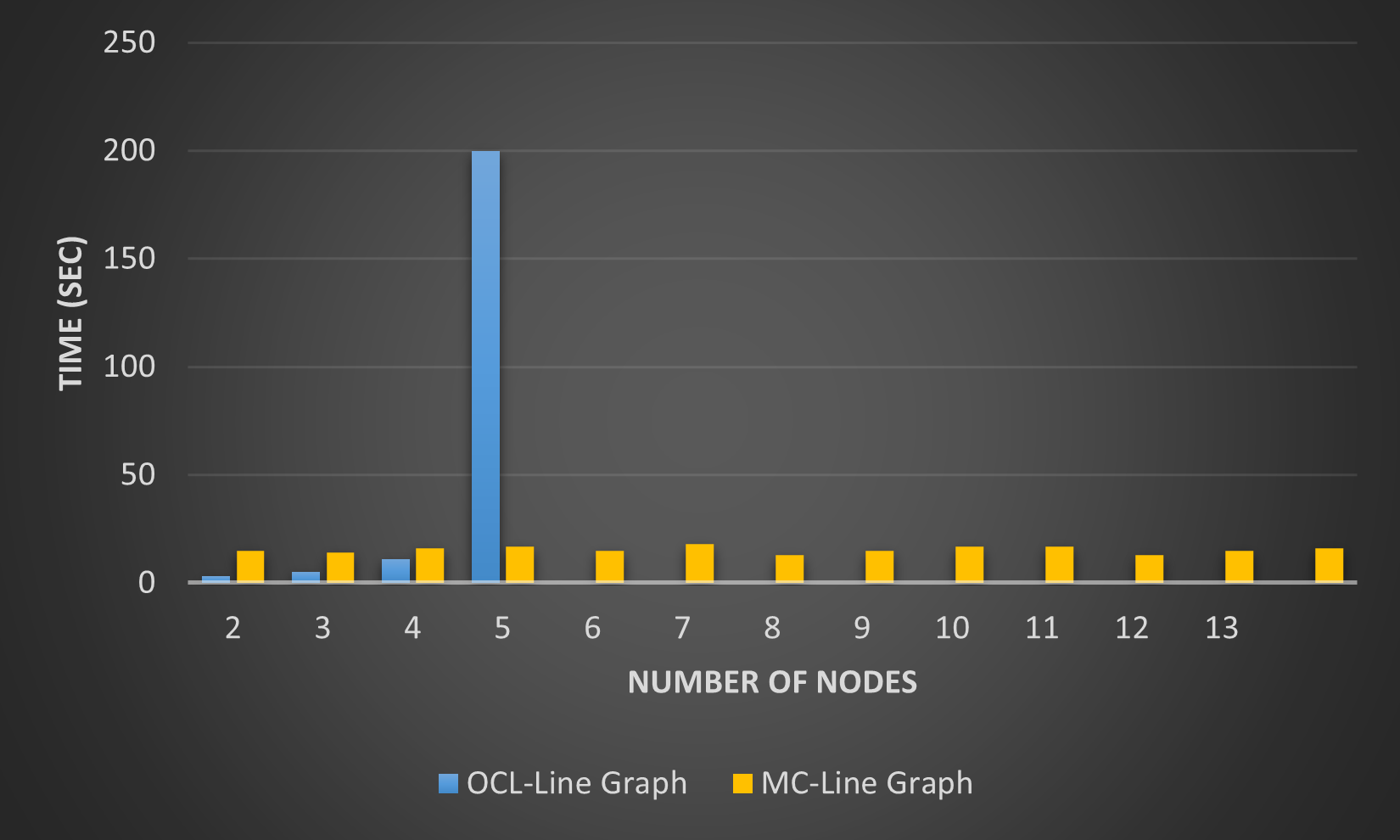}}
\caption{(a) time growth of finding a solution with respect to number of agents for complete graph for OpenOCL tool (OCL) and MonteCarlo (MC).(b)time growth of finding a solution with respect to number of agents for cycle graph graph for OpenOCL tool (OCL) and MonteCarlo (MC). (c) time growth of finding a solution with respect to number of agents for line graph graph for OpenOCL tool (OCL) and MonteCarlo (MC).}
    \label{t_compl}
\end{figure*}

In this last section, we would like to report on some observations on the time complexity of the two methods, OCL Strategy and Threshold Strategy, that were used in subsection~\ref{se:simulation}. We study the time growth, as a function of the number of agents in a network, of the two algorithms that computes the solution. We run the two algorithms OpenOCL (OCL) and MonteCarlo (MC) for three different types of graphs, complete, cycle, and line graphs as can be seen in Figure~\ref{t_compl}. Our observations from these data can tell us that the time complexity of OCL can be worse than $\omega(N^{\log(N))}$, where $\omega({\cdot})$ represents the best case scenario in the time complexity theory. On the other hand, time complexity of using MonteCarlo method in its worst case is $o(\ell)$ where $\ell$ represents number of iterations and $o(\cdot)$ defines the worst case scenario in time complexity theory. 

\begin{figure}[h]
	\centering
	\includegraphics[width=0.9\linewidth]{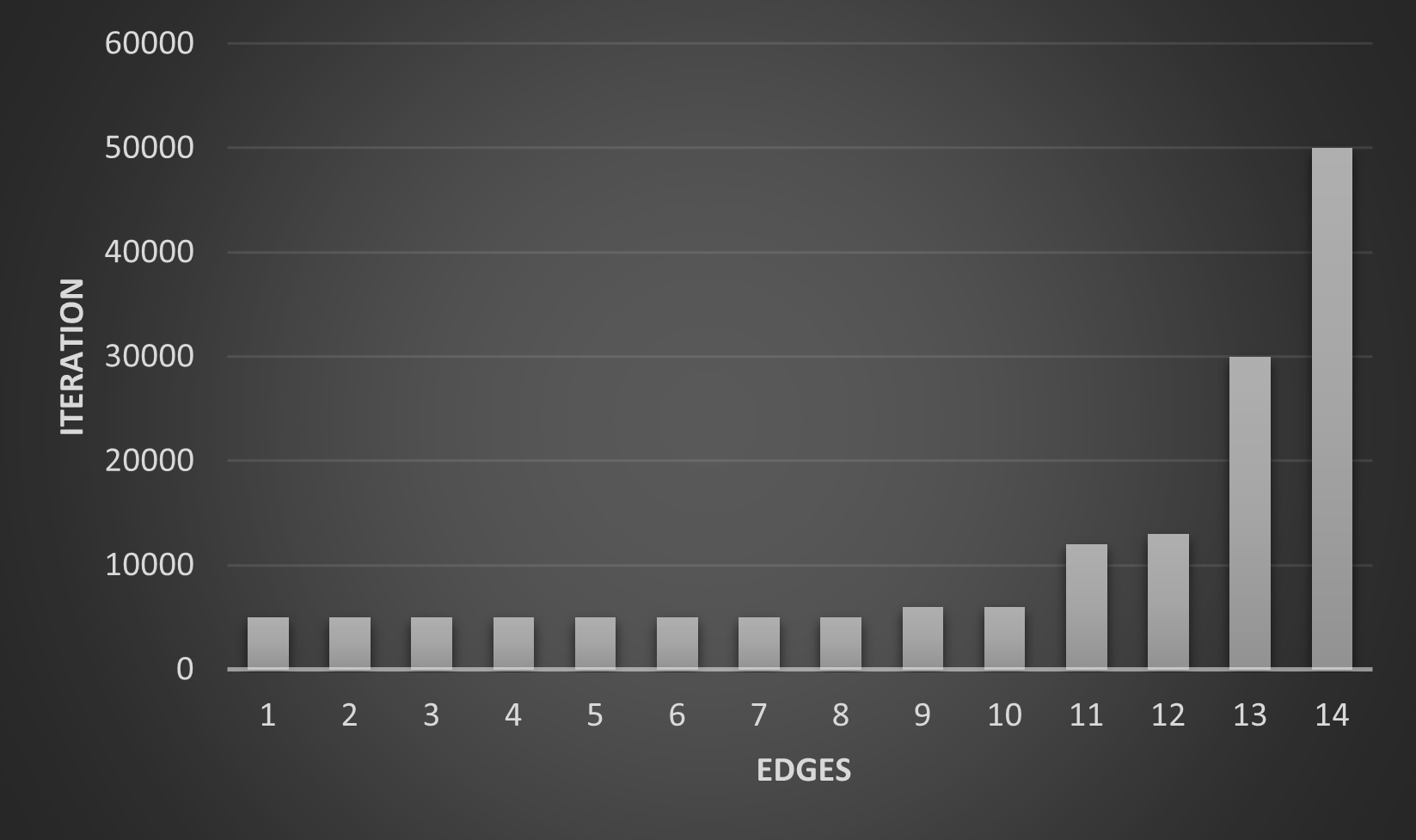}
	\caption{Data collected for the relation between number of iterations needed when using MC method to find the solution vs number of edges being controlled. }
\label{fig:edge_it}
\end{figure}

We also observed that when using MonteCarlo method, there is a relation between the number of edges that are being controlled and the number of iterations needed to find the solution, as can be seen in Figure~\ref{fig:edge_it}. 

\section{\textbf{Conclusions}}
In this paper, we formulate a real world problem, of how to optimally control  individual interactions with their neighbors. 
We formulated the problem in two different scenarios, first we formalize a distributed optimal control. However, due to combinatorial complexity of the problem finding optimal policies are intractable even for small network. We present a bottom-up approach to solve the problem, where a local control strategy for each agent is formalized. Interestingly, under some valid assumptions, we show that the control strategy to the problem to be a form of optimal infection threshold of each neighbor. To compare our solution against some other solutions, we solved the relaxed problem numerically using Open OCL software package, 
Finally, we use well-mixed homogeneous population case, and we showed that this assumption is as good even though it simplify the analysis greatly. 

\section*{Acknowledgement}
This material is based upon research supported in part by the National Science Foundation under grant number 2028523 and the U.S. Office of Naval Research under award number N00014-22-1-2207.

{\scriptsize
  \bibliographystyle{ieeetr}
  \bibliography{Mohammad_Main}
%
}

\section{Appendix}
\subsection{Proof of Theorem \eqref{th:solut1}}
Now, setting up an optimal control problem as in Table 3.2-1 of~\cite{FllDlvVls2012}, we have the cost function 
\begin{align*}
J = \int_0^T \sum_{j \in \NN_1} c_{1j}u_{1j}(t) dt,
\end{align*}
and the terminal constraint is 
\begin{align*}
\psi(x(T)) &= \x^S_1(T) - \threshold_1 = 0 .
\end{align*}
This means the Hamiltonian is
\begin{align*}
H &= \sum_{j \in \NN_1} c_{1j}u_{1j} + \sum_{i \in \mathcal{V}}\lambda_{\x_i^S}\dot{\x}_i^S+\sum_{i \in \mathcal{V}}\lambda_{\x_i^I} \dot{\x}_i^I
\end{align*}
where $\lambda \in \R^{2N}$ is the costate, $\lambda_{\x_i^S} \in \{\lambda_1,\dots,\lambda_{N}\}$ and $\lambda_{\x_i^I} \in \{\lambda_{N+1},\dots,\lambda_{2N}\}$. 
This means that the costate dynamics is
\begin{align*}
\dot{\lambda}_{\x_i^S} &= -\frac{\partial H}{\partial \x_i^S}\\
\dot{\lambda}_{\x_i^I} &= -\frac{\partial H}{\partial \x_i^I}^T
%
\end{align*}

The other boundary condition is
\begin{align*}
v\frac{\partial \psi}{\partial x}^T(T) - \lambda(T) =& 0,\\
\end{align*}
which implies $\lambda_1(T) = v$, for some $v \in \R$, and $\lambda_i(T)=0$, for $i = 2,\dots,2N$.

We need to use Pontryagin's minimum principle (PMP) (Chapter 5.2 of~\cite{FllDlvVls2012}) to get an expression for the optimal control $u^*$
\begin{align}
H(\x^{S*},\x^{I*},u^*,\lambda^*) &\leq H(\x^{S*},\x^{I*},u,\lambda^*)\nonumber\\
\sum_{j \in \NN_1} c_{1j}u^*_{1j} + \sum_{i \in \mathcal{V}}\lambda_{\x_i^{S*}}\dot{\x}_i^{S*}+\lambda_{\x_i^{I*}}\dot{\x}_i^{I*} &\leq\nonumber\\ 
\sum_{j \in \NN_1} c_{1j}u_{1j} + \sum_{i \in \mathcal{V}}\lambda_{\x_i^{S*}}\dot{\x}_i^{S*}+\lambda_{\x_i^{I*}}\dot{\x}_i^{I*},
\end{align}
where $\x^{S*}=[\x_1^{S*},\dots,\x_N^{S*}]^T, \, \x^{I*}=[\x_1^{I*},\dots,\x_N^{I*}]^T$. After some analysis and simplifications, the inequality becomes as follows,
\begin{align}\label{eq:pontryaginGeneralCase2}
\sum_{j\in \NN_1}(c_{1j}+\sum_{i \in \mathcal{V}}(\lambda_{\x_i^{S*}}-\lambda_{\x_i^{I*}}) \x_i^{S*} \x_j^{I*}) u_{1j}^* &\leq  \nonumber\\
\sum_{j\in \NN_1}(c_{1j}+\sum_{i \in \mathcal{V}}(\lambda_{\x_i^{S*}}-\lambda_{\x_i^{I*}}) \x_i^{S*}\x_j^{I*} )u_{1j},
\end{align}
for all $j\in\NN_1$, and using the fact that the $u_{1j}$'s are independent. Therefore, if $c_{1j}+\sum_{i \in \mathcal{V}}(\lambda_{\x_i^{S*}}-\lambda_{\x_i^{I*}}) \x_i^{S*} \x_j^{I*} < 0$, then $u_{1j}^*=1$, and if $c_{1j}+\sum_{i \in \mathcal{V}}(\lambda_{\x_i^{S*}}-\lambda_{\x_i^{I*}}) \x_i^{S*} \x_j^{I*} > 0$, then $u_{1j}^*=0$. Therefore, we can write the input as
\begin{align}\label{inputWswitchingFunc2}
u_{1j}^* &= \frac{1}{2} - \frac{1}{2}\text{sgn }\left( c_{1j}+\sum_{i \in \mathcal{V}}(\lambda_{\x_i^{S*}}-\lambda_{\x_i^{I*}}) \x_i^{S*} \x_j^{I*}\right).
\end{align}

\subsection{Proof of Theorem \eqref{th:solut}}
For convenience, we define $\chi^I = [\x_2^I,\dots,\x_N^I]^T \in \R^{N-1}$, $\chi^S = [\x_2^S,\dots,\x_N^S]^T \in \R^{N-1}$, and $\chi \triangleq [\chi^{IT}, \chi^{ST}]^T$. Setting up an optimal control problem as in Table 3.2-1 of~\cite{FllDlvVls2012},we have the cost function 
\begin{align*}
J = \int_0^T \sum_{j \in \NN_1} c_{1j}u_{1j}(t) dt,
\end{align*}
and the terminal constraint is 
\begin{align*}
\psi(x(T)) &= \x^S_1(T) - \threshold_1 = 0 .
\end{align*}
This means the Hamiltonian is
\begin{align*}
H &= \sum_{j \in \NN_1} c_{1j}u_{1j} + \lambda_1\dot{\x}_1^S+\lambda_\chi^T \dot{\chi}
\end{align*}
where $\lambda \in \R^{2N-1}$ is the costate, and $\lambda_\chi \triangleq [\lambda_2,\dots,\lambda_{2N-1}]^T$. 
This means that the costate dynamics is
\begin{align*}
\dot{\lambda}_1 &= -\frac{\partial H}{\partial \x_1^S}\\
\dot{\lambda_\chi} &= -\frac{\partial H}{\partial \chi}^T
%
\end{align*}

The other boundary condition is
\begin{align*}
v\frac{\partial \psi}{\partial x}^T(T) - \lambda(T) =& 0,\\
\end{align*}
which implies $\lambda_1(T) = v$, for some $v \in \R$, and $\lambda_i(T)=0$, for $i = 2,\dots,2N-1$.

We need to use PMP (Chapter 5.2 of~\cite{FllDlvVls2012}) to get an expression for the optimal control $u^*$
\begin{align}
H([\x_1^{S*T},\chi^{*T}]^T,u^*,\lambda^*) &\leq H([\x_1^{S*T},\chi^{*T}]^T,u,\lambda^*)\nonumber\\
\sum_{j \in \NN_1} c_{1j}u_{1j}^* + \lambda_1^* \dot{\x}_1^{S*}+\lambda_\chi^{T*} \dot{\chi}^* &\leq \sum_{j \in \NN_1} c_{1j}u_{1j}+ \lambda_1^* \dot{\x}_1^{S*}+\lambda_\chi^{T*} \dot{\chi}^*,\nonumber
\end{align}
since the only term of the states that depends on the input is $\dot{\x}_1^{S*}$, the inequality becomes as follows,
\begin{align}\label{eq:pontryaginGeneralCase}
\sum_{j\in \NN_1}(c_{1j}+\lambda_1^* \beta \x_1^{S*}\x_j^{I*}) u_{1j}^* &\leq  \sum_{j\in \NN_1}(c_{1j}+\lambda_1^* \beta \x_1^{S*}\x_j^{I*} )u_{1j}
\end{align}
for all $j\in\NN_1$, and using the fact that the $u_{1j}$'s are independent. Therefore, if $c_{1j}+\lambda_1^* \beta \x_1^{S*}\x_j^{I*} < 0$, then $u_{1j}^*=1$, and if $c_{1j}+\lambda_1^* \beta \x_1^{S*}\x_j^{I*} > 0$, then $u_{1j}^*=0$. Therefore, we can write the input as
\begin{align}\label{inputWswitchingFunc}
u_{1j} &= \frac{1}{2} - \frac{1}{2}\text{sgn }\left( c_{1j}+\lambda_1 \beta \x_1^{S}\x_j^{I}\right).
\end{align}
Now, note that
\begin{align*}
\dot{\lambda}_1 &= -\frac{\partial H}{\partial \x_1^S} = -\lambda_1 \frac{\partial}{\partial \x_1^S} \dot{\x}_1^S,
\end{align*}
because only $\dot{\x}_1^S$ depends on $\x_1^S$. Now, from~\eqref{eq:dynamics1}, we have $\frac{\partial}{\partial \x_1^S} \dot{\x}_1^S = -\sum_{j\in\NN_1}{\beta}(1-u_{1j}(t))\x_j^I$, which means
\begin{align*}
\dot{\lambda}_1 &= -\frac{\lambda_1\dot{\x}_1^S}{\x_1^S}.
\end{align*}
This implies that $0 = \dot{\lambda}_1\x_1^S + \lambda_1\dot{\x}_1^S =\frac{d}{dt}(\lambda_1(t)\x_1^S(t))$, indicating that the quantity $\lambda_1(t)\x_1^S(t)$ is constant with respect to time. 
At the terminal time we know, $\lambda_1(T)\x_1^S(T)=v\threshold_1$, and so $\lambda_1(t)\x_1^S(t)=v\threshold_1$, $\forall t \in [0,T]$. Therefore, the optimal input~\eqref{inputWswitchingFunc} can be written as 
\begin{align}
u_{1j}&=\begin{cases}
1,\quad &\text{if}~ v\x^{I}_j< \frac{-c_{1j}}{\beta \threshold_1}, \nonumber\\ 
0,\quad &\text{if}~ v\x^{I}_j> \frac{-c_{1j}}{\beta \threshold_1}.
\end{cases}
\end{align}
Note that $c_{1j},\beta,\threshold_1,v$ are all constant in time and $c_{1j},\beta,\threshold_1,\x^{I}_j>0$, so we have the trivial solution $u_{1j}^*(t)=0, \forall t \in [0,T]$ if $v \geq 0$. In the nontrivial case where $v<0$, then, we can write the optimal solution as
\begin{align}
u_{1j}&=\begin{cases}
1,\quad &\text{if}~ \x^{I}_j> \mathcal{T}_j^I, \nonumber\\
0,\quad &\text{if}~ \x^{I}_j< \mathcal{T}_j^I,
\end{cases}
\end{align}
where $\mathcal{T}_j^I \triangleq \frac{-c_{1j}}{\beta \threshold_1 v}$. Noting that, since $\x^{I}_j \in [0,1]$, we can recover the trivial case with a proper selection of $\mathcal{T}_j^I$, we have proven that the optimal input has the threshold form.

\subsection{Proof of Theorem \eqref{th:solution2}}
Setting up an optimal control problem as in Table 3.2-1 of~\cite{FllDlvVls2012}, the cost function with lumped population becomes  
\begin{align*}
J = \int_0^T c_{av}u(t) dt,
\end{align*}
where $c_{av}=\frac{N-1}{|\NN_1|}\sum_{j\in \NN_1}c_{1j}$, is the average of all edge's associated cost.

The terminal constraint is 
\begin{align*}
\psi(x(T)) &= \x^S_1(T) - \threshold_1 = 0 .
\end{align*}
This means the Hamiltonian is
\begin{align*}
H &=  c_{av} u + \lambda_1\dot{\x}_1^S+\lambda_2 \dot{I_{LP}}+ \lambda_3 \dot{S_{LP}}+
\end{align*}
where $\lambda \in \R^{3}$ is the costate, and $\lambda \triangleq [\lambda_1,\lambda_2,\lambda_3]^T$. 
This means that the costate dynamics is
\begin{align*}
\dot{\lambda}_1 &= -\frac{\partial H}{\partial \x_1^S}\\
\dot{\lambda_2} &= -\frac{\partial H}{\partial I_{LN}}^T\\
\dot{\lambda_3} &= -\frac{\partial H}{\partial S_{LN}}^T.
%
\end{align*}

The other boundary condition is
\begin{align*}
v\frac{\partial \psi}{\partial x}^T(T) - \lambda(T) =& 0,\\
\end{align*}
which implies $\lambda_1(T) = v$, for some $v \in \R$, and $\lambda_i(T)=0$, for $i = 2,3 $.

We need to use PMP (Chapter 5.2 of~\cite{FllDlvVls2012}) to get an expression for the optimal control $u^*$
\begin{align}
H(\x_1^{S*},I^*_{LN},S^*_{LN},u^*,\lambda^*) &\leq H(\x_1^{S*},I^*_{LN},S^*_{LN},u,\lambda^*)\nonumber\\
c_{av}u^* + \lambda_1^* \dot{\x}_1^{S*}+\lambda_2^{*} \dot{I}_{LN}^*+\lambda_3^{*} \dot{S}_{LN}^* &\leq \nonumber\\ c_{av}u + \lambda_1^* \dot{\x}_1^{S*}+\lambda_2^{*} \dot{I}_{LN}^*+\lambda_3^{*} \dot{S}_{LN}^*,\nonumber
\end{align}
since the only term of the states that depends on the input is $\dot{\x}_1^{S*}$, the inequality becomes as follows,
\begin{align}\label{eq:pontryaginGeneralCase3}
(c_{av}+\lambda_1^* \beta \x_1^{S*}{I^*_{LN}}) u^* &\leq  (c_{av}+\lambda_1^* \beta \x_1^{S*}{I^*_{LN}} )u.
\end{align}
Now, if $c_{av}+\lambda_1^* \beta \x_1^{S*}{I^*_{LN}}< 0$, then $u^*=1$, and if $c_{av}+\lambda_1^* \beta \x_1^{S*}{I^*_{LN}} > 0$, then $u^*=0$. Therefore, we write the input as
\begin{align}\label{inputWswitchingFunc3}
u^*&= \frac{1}{2} - \frac{1}{2}\text{sgn }\left(c_{av}+\lambda_1^* \beta \x_1^{S*}I^*_{LN}\right).
\end{align}
Now, note that
\begin{align*}
\dot{\lambda}_1 &= -\frac{\partial H}{\partial \x_1^S} = -\lambda_1 \frac{\partial}{\partial \x_1^S} \dot{\x}_1^S,
\end{align*}
because only $\dot{\x}_1^S$ depends on $\x_1^S$. Now, from~\eqref{eq:dynamics1}, we have $\frac{\partial}{\partial \x_1^S} \dot{\x}_1^S = -{\beta}(1-u(t))I_{LN}$, which means
\begin{align*}
\dot{\lambda}_1 &= -\frac{\lambda_1\dot{\x}_1^S}{\x_1^S}.
\end{align*}
This implies that $0 = \dot{\lambda}_1\x_1^S + \lambda_1\dot{\x}_1^S =\frac{d}{dt}(\lambda_1(t)\x_1^S(t))$, indicating that the quantity $\lambda_1(t)\x_1^S(t)$ is constant with respect to time. 
At the terminal time we know, $\lambda_1(T)\x_1^S(T)=v\threshold_1$, and so $\lambda_1(t)\x_1^S(t)=v\threshold_1$, $\forall t \in [0,T]$. Therefore, the optimal input~\eqref{inputWswitchingFunc3} can be written as 
\begin{align}
u^*&=\begin{cases}
1,\quad &\text{if}~ vI_{LN}< \frac{-c_{av}}{\beta \threshold_1}, \nonumber\\ 
0,\quad &\text{if}~ vI_{LN}> \frac{-c_{av}}{\beta \threshold_1}.
\end{cases}
\end{align}
Note that $c_{av},\beta,\threshold_1,v$ are all constant in time and $c_{av},\beta,\threshold_1,I_{LN}>0$, so we have the trivial solution $u^*(t)=0, \forall t \in [0,T]$ if $v \geq 0$. In the nontrivial case where $v<0$, then, we can write the optimal solution as
\begin{align}
u^*&=\begin{cases}
1,\quad &\text{if}~ I_{LN}> I^*, \nonumber\\
0,\quad &\text{if}~ I_{LN}< I^*,
\end{cases}
\end{align}
where $I^* \triangleq \frac{-c_{av}}{\beta \threshold_1 v}$. Noting that, since $I^* \in [0,1]$, we can recover the trivial case with a proper selection of $I^*$, we have proven that the optimal input has the threshold form.

\begin{IEEEbiography}[{\includegraphics[width=1in,height=1.25in,clip,keepaspectratio]{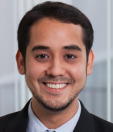}}]{Cameron Nowzari}
received the Ph.D. in Mechanical Engineering from the University of California, San Diego in September 2013. He then held a postdoctoral position with the Electrical and Systems Engineering Department at the University of Pennsylvania until 2016. He is currently an Assistant Professor with the Electrical and Computer Engineering Department at George Mason University, in Fairfax, Virginia. He has received several awards including the American Automatic Control Council’s O. Hugo Schuck Best Paper Award, the IEEE Control Systems Magazine Outstanding Paper Award, and the International Conference on Data Mining Best Paper Award. His current research interests include dynamical systems and control, distributed coordination algorithms, robotics, event- and self-triggered control, Markov processes, network science, spreading processes on networks, and the Internet of Things.
\end{IEEEbiography}
\begin{IEEEbiography}[{\includegraphics[width=1in,height=1.25in,clip,keepaspectratio]{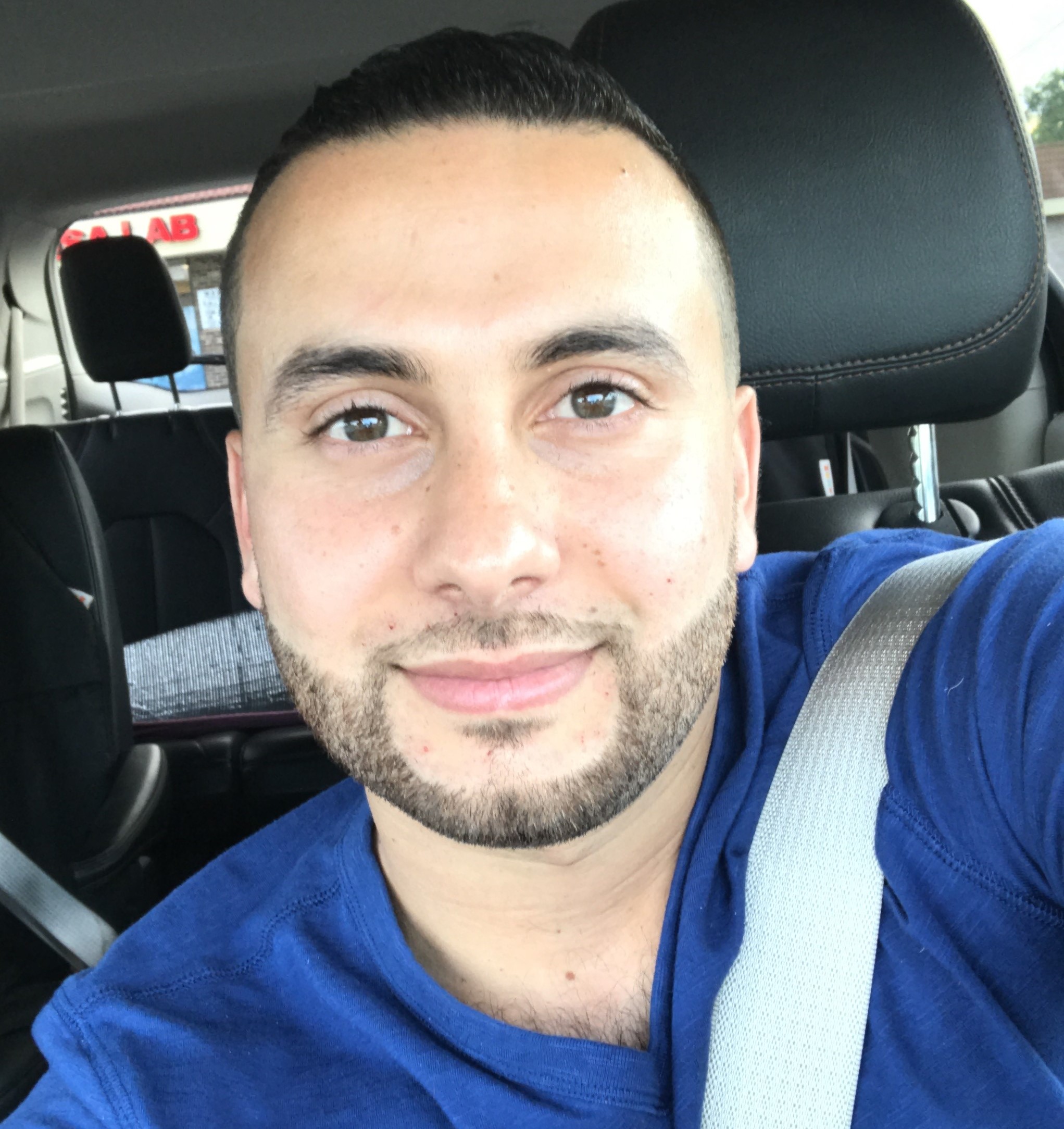}}]{Mohammad Mubarak}
received the BS in Electrical Engineering from Al-Balqa Applied University in 2014. He then worked as a research and development engineer in Reyadah Electronics until 2017. He is currently pursuing a PhD in Electrical and Computer Engineering from George Mason University. His research interest including multi-agent systems, distributed control, dynamical systems control, spreading process on network, and Markov Process.
\end{IEEEbiography}

\end{document}